\documentclass[onecolumn,final,a4paper]{elsarticle}
\graphicspath{{Figures/}}

\usepackage[utf8]{inputenc}
\usepackage{amsmath,amssymb,amsfonts,latexsym,stmaryrd}
\usepackage{algorithmic}
\usepackage{booktabs}

\usepackage[colorlinks=true,breaklinks=true,linkcolor=lightblue,citecolor=lightblue]{hyperref}
\usepackage{stackengine}
\usepackage{mathrsfs}
\usepackage{rotating} 
\usepackage{amssymb,amsmath,bbm,amsfonts,latexsym,subfigure}
\usepackage{amsthm}
\usepackage{graphicx,float,epsfig,color,fancyhdr}
\usepackage{amsopn}

\usepackage[utf8]{inputenc}
\usepackage{graphicx}
\usepackage{wrapfig}

\usepackage{amsfonts}
\usepackage{amsmath}
\numberwithin{equation}{section}
\usepackage{mathtools}

\usepackage[ruled,vlined,linesnumbered]{algorithm2e}
\let\oldnl\nl
\newcommand{\nonl}{\renewcommand{\nl}{\let\nl\oldnl}}%

\usepackage[utf8]{inputenc}
\usepackage{graphicx}
\usepackage{wrapfig}

\usepackage{amsfonts}
\usepackage{amsmath}
\numberwithin{equation}{section}
\usepackage{mathtools}

\usepackage{makecell}

\usepackage[ruled,vlined,linesnumbered]{algorithm2e}
\let\oldnl\nl

\newcommand{\tenprod}{\circ}
\newcommand{\ten}[1]{\mathcal{#1}}
\newcommand{\mat}[1]{\mathbf{#1}}

\newcommand{\PGRAPH}[1]{\noindent\textbf{#1}}

\newcommand{\HAT}[1]{\widetilde{}}

\usepackage{lipsum}

\newtheorem{lemma}{Lemma}[section]

\newtheorem{theorem}{Theorem}[section]

\def\bA{\mathbf{A}}
\def\bM{\mathbf{M}}

\def\bH{\mathbf{H}}

\def\bB{\mathbf{B}}

\def\bD{\mathbf{D}}
\def\bE{\mathbf{E}}
\def\bF{\mathbf{F}}

\def\bZero{\mathbf{0}}

\def\bv{\mathbf{v}}

\def\bff{\mathbf{f}}

\def\bx{\mathbf{x}}

\def\bI{\boldsymbol{I}}

\def\bD{{\bf D}}

\def\bV{\boldsymbol{\mathcal{V}}}
\def\bScal{\boldsymbol{\mathcal{S}}}

\def\bS{\mathbf{S}}

\def\curl{\mathop{\mathbf{curl}}\nolimits}

\def\div{\mathop{\mathrm{div}}\nolimits}

\def\J{\mathbf{J}}

\newcommand{\vertiii}[1]{{\left\vert\kern-0.25ex\left\vert\kern-0.25ex\left\vert #1 
    \right\vert\kern-0.25ex\right\vert\kern-0.25ex\right\vert}}

\newcommand{\TTf}{\text{TT}}

\begin{document}

\begin{frontmatter}

  \title{ Space-Time Spectral Collocation Tensor Network Approach for Maxwell's Equations}

  \author[TDIV,IIT-KGP]{Dibyendu Adak}
  \author[TDIV]{Rujeko Chinomona}
  \author[TDIV]{Duc P. Truong}
  \author[TDIV]{Oleg Korobkin}
  \author[TDIV]{Kim {\O}. Rasmussen}
  \author[TDIV]{Boian S. Alexandrov}

  \address[TDIV]{Theoretical Division,
    Los Alamos National Laboratory, Los Alamos, NM 87545, USA
  }
  \address[IIT-KGP]{Indian Institute of Technology-Kharagpur, WB, India}
  

\begin{abstract}
In this work, we develop a space--time Chebyshev spectral collocation method for three-dimensional Maxwell's
equations and combine it with tensor-network techniques in Tensor-Train (TT) format. 
Under constant material parameters, the Maxwell system is reduced to a vector wave equation for the electric field, which we discretize globally in space and time using a staggered spectral collocation scheme. 
The staggered polynomial spaces are designed so that the discrete curl and divergence operators preserve the divergence-free constraint on the magnetic field. 
The magnetic field is then recovered in a space--time post-processing step via a discrete version of Faraday's law. 
The global space–time formulation yields a large but highly structured linear system, which we approximate in low-rank TT-format directly from the operator and data, without assuming that the forcing is separable in space and time.
We derive condition-number bounds for the resulting operator and prove spectral convergence estimates for both the electric and magnetic fields.
Numerical experiments for three-dimensional electromagnetic test problems confirm the theoretical convergence rates and show that the TT-based solver maintains accuracy with approximately linear complexity in the number of grid points in space and time.

\end{abstract}

\end{frontmatter}

\section{Introduction}

The accurate and efficient numerical solution of Maxwell's equations remains a cornerstone challenge in computational electromagnetics, with applications spanning satellite communications, medical imaging, radar systems, telecommunication devices, and numerous other technologies critical to modern society. Maxwell's equations describe the propagation and scattering of electromagnetic waves, and their time-dependent nature requires careful treatment to preserve physical properties while maintaining computational efficiency. As electromagnetic devices and systems become increasingly complex and operate at higher frequencies, the demand for high-fidelity simulations that can resolve fine-scale features across large computational domains has grown substantially.

Many numerical approaches for time-dependent partial differential equations combine spectral methods for spatial discretization with low-order finite difference schemes for temporal discretization \cite{hussaini1989spectral}. While spectral methods offer exponential convergence for smooth solutions in the spatial domain, this advantage is undermined when temporal errors dominate, reducing the overall convergence to first or second order. To address this limitation, space-time spectral collocation methods have been developed \cite{lui2020chebyshev, lui2017legendre, lui2021spectral}, where both spatial and temporal dimensions are discretized using high-order spectral approximations. These methods achieve exponential convergence globally but require solving for all time steps simultaneously, substantially increasing the dimensionality of the problem and introducing significant computational challenges.

For Maxwell's equations specifically, preserving the divergence-free constraint of the magnetic field ($\nabla \cdot \mathbf{B} = 0$) poses a fundamental numerical challenge. Classical approaches have addressed this through various means. Yee's pioneering finite-difference time-domain scheme \cite{yee1966numerical} introduced spatial staggering of field components on a Cartesian grid, a principle that has been extended to fourth-order accurate compact schemes \cite{yefet2000fourth, yefet2001staggered}. Discontinuous Galerkin methods \cite{araujo2017stability, camargo2020hdg} and multidomain methods \cite{collino1997new} have also been developed to maintain divergence-free properties while allowing for complex geometries and material interfaces. In the spectral element context, Monk \cite{monk1992analysis} analyzed finite element discretizations of Maxwell's equations, while more recent work has explored Legendre-tau Chebyshev collocation spectral element methods with Runge-Kutta time integration \cite{niu2023legendre} and divergence-free finite element schemes for static problems \cite{huang2012divergence}. Despite these advances, the computational cost of these methods remains substantial, particularly for three-dimensional problems requiring fine resolution over long time intervals.

A more fundamental obstacle confronts all traditional approaches to solving Maxwell's equations in three spatial dimensions plus time: the curse of dimensionality \cite{bellman1966dynamic}. The number of degrees of freedom grows exponentially with dimension, rendering fine-resolution simulations prohibitively expensive or entirely infeasible. Consider, for example, full waveform inversion problems in geophysics or radar imaging, where the number of grid points per time step can reach $M = 10^9$ or more \cite{krebs2009fast}. With $N$ time steps required to simulate the wave propagation, one must store and manipulate $MN$ floating point numbers, quickly exceeding available memory even on modern high-performance computing systems. This exponential scaling defeats traditional computational approaches and forces difficult trade-offs. One option is to employ reduced-order models (ROMs) \cite{gunzburger2007reduced, leugering2014trends, danis2025tensor}, which construct low-dimensional approximations of the solution manifold. While ROMs excel for certain problem classes, they often require significant domain expertise, are tailored to specific applications, and demand invasive modifications to existing simulation codes. Another option is checkpointing, where only selected time steps are stored and intermediate states are recomputed as needed, but this approach multiplies the computational cost and still faces memory limitations for large-scale problems.

Tensor network methods have emerged as a promising alternative for mitigating the curse of dimensionality in high-dimensional problems. These methods exploit the observation that many high-dimensional functions of practical interest, despite having exponentially many degrees of freedom in principle, can be well approximated by low-rank tensor decompositions. Among tensor formats, the tensor-train (TT) decomposition introduced by Oseledets \cite{oseledets2011tensor} has gained widespread adoption because it can accurately represent complex functions while remaining computationally feasible. The TT format represents a $d$-dimensional tensor as a product of three-dimensional cores, effectively implementing a discrete separation of variables \cite{bachmayr2016tensor}. Crucially, when the underlying function admits a low-rank structure, the TT representation achieves linear complexity in both storage and computation with respect to dimension and discretization level. This strongly contrasts the exponential scaling of conventional methods.

The utility of tensor-train methods rests on efficient algorithms that can construct and manipulate TT representations without ever forming the full tensor. The TT-cross interpolation algorithm \cite{oseledets2010tt, savostyanov2011fast, sozykin2022ttopt} addresses the construction problem by adaptively selecting tensor fibers that capture the essential structure, building the decomposition from carefully chosen one-dimensional slices rather than the complete high-dimensional array. For solving linear systems in TT format, alternating optimization methods \cite{oseledets2012solution, dolgov2014alternating, holtz2012alternating} have proven effective, iteratively updating each core while fixing the others. These algorithms have been successfully applied to quantum many-body problems, chemical master equations, and various PDE systems. Moreover, the Quantized Tensor-Train format \cite{kazeev2013multilevel} provides additional compression for operators with special structure, such as the banded matrices arising from spectral discretizations.

Recent work has begun to explore the combination of tensor network methods with high-order discretizations for PDEs. Applications include the Boltzmann neutron transport equation \cite{truong2023tensor}, where the TT format enables simulations in the full six-dimensional phase space; the mimetic finite difference method for Maxwell's equations \cite{manzini2023tensor}; and space-time spectral collocation for nonlinear convection-diffusion equations \cite{adak2024tensor, adak2024tensorA}, where the TT approach has demonstrated dramatic computational savings while preserving spectral accuracy. These successes suggest that tensor methods may offer a path forward for large-scale electromagnetic simulations that would be intractable with conventional approaches.

In this work, we develop a novel space-time spectral collocation method for Maxwell's equations that leverages tensor-train decomposition to overcome the curse of dimensionality while maintaining exponential convergence. Our approach begins with a global space-time formulation where both spatial and temporal dimensions are discretized using Chebyshev-Gauss-Lobatto collocation points. This yields a large block-structured linear system that encapsulates the entire time evolution, making all time steps available simultaneously for tensor compression. Under the assumption of constant material parameters, we reformulate Maxwell's system as a vector wave equation for the electric field \cite{monk1992analysis}, allowing us to solve for each electric field component independently. The magnetic field is then recovered through post-processing using a discretization of Faraday's law. To preserve the divergence-free property of the magnetic field without explicit constraint enforcement, we use a spectral analogue of classical staggered grids: electric and magnetic field components live in different polynomial spaces with carefully chosen degrees of freedom in each coordinate direction. This spectral staggering ensures that the discrete curl operator maps between compatible spaces and that the discrete identity $\nabla \cdot (\nabla \times \cdot) \approx 0$ holds to spectral accuracy.

The Kronecker product structure of the discrete operators arising from this space-time formulation proves naturally compatible with tensor-train representations. Each spatial and temporal derivative operator can be expressed as a Kronecker product of one-dimensional matrices, and the resulting discrete wave equation inherits this structure. We construct TT representations of the solution, operators, and boundary data using TT-cross interpolation, solving the compressed system using alternating minimal energy methods \cite{dolgov2014alternating}. The theoretical analysis establishes rigorous foundations for the method, including condition number estimates showing $\mathcal{O}(N^4)$ scaling for the wave operator and $\mathcal{O}(N^2)$ scaling for the temporal operator, as well as proofs of exponential convergence for both electric and magnetic fields. We also prove that both Gauss's law for the electric field and the divergence-free constraint for the magnetic field are satisfied with exponential accuracy.

Numerical experiments validate these theoretical results and demonstrate the practical efficiency of the method. Across test cases with manufactured solutions of varying tensor rank, we consistently observe exponential convergence matching the full-grid spectral method while achieving computational speedups of $10^7$ to $10^8$ for problems with approximately $2 \times 10^5$ space-time degrees of freedom. Critically, the tensor-train method extends well beyond the memory limitations of full-grid computation, successfully handling discretizations that would be completely intractable with conventional approaches.


Our paper is organized as follows.
Section~\ref{sec:math:model:num:dis} introduces the mathematical model, derives the second-order
formulation for the electric field, and presents the space--time Chebyshev spectral collocation
discretization with staggered spaces. Section~\ref{sec:wave_equation} develops the matrix
formulation of the wave equation and describes the magnetic field reconstruction via
post-processing. Section~\ref{sec:conditioning_convergence} establishes condition number
estimates and spectral convergence rates for both electric and magnetic fields.
Section~\ref{sec:Tensor-Networks} reviews tensor network concepts and describes the tensorization
procedure for both the electric field solution and magnetic field recovery.
Section~\ref{sec:numerical_experiments} presents numerical experiments demonstrating exponential
convergence and computational speedups compared to full-grid methods.
Section~\ref{sec:conclusions} offers concluding remarks and discusses future directions.
\ref{APP:A_map_tt} provides technical details on the construction of mapping operators
and boundary tensors.

\section{Mathematical model and discretization}
\label{sec:math:model:num:dis}

Let $\Omega \subset \mathbb{R}^3$ be a cubic domain with boundary $\partial \Omega$. We consider
Maxwell's equations on $\Omega$ for the electric field
$\bE = (E_x,E_y,E_z)$ and the magnetic flux density
$\bB = (B_x,B_y,B_z)$:
\begin{align}
       & \frac{\partial \bB}{\partial t} = -\nabla \times \bE
         \qquad\qquad\qquad\quad \text{in } \Omega,\ \forall t \in (0,T], \label{model:original:1}\\
       & \frac{\partial \bE}{\partial t} = c^2 \nabla \times \bB -\frac{1}{\epsilon_0} \J
         \qquad\qquad \text{in } \Omega,\ \forall t \in (0,T], \label{model:original:2} \\
       & \nabla \cdot \bE = \rho
         \qquad\qquad\qquad\qquad\quad\ \text{in } \Omega, \label{model:original:3}\\
       & \nabla \cdot \bB = 0
         \qquad\qquad\qquad\qquad\quad\ \text{in } \Omega,  \label{model:original:4}
\end{align}
where $\J$ is the current density, $\rho$ the charge density, $\epsilon_0$ the dielectric constant, and $c$ the speed of light. We consider the model problem
\eqref{model:original:1}–\eqref{model:original:4} with initial conditions
$\bE(0,\cdot) = \bE^0$, $\bB(0,\cdot) = \bB^0$ and inhomogeneous boundary conditions for both
the electric field $\bE$ and the magnetic field $\bB$. The implementation of inhomogeneous
boundary conditions for the wave equation is not straightforward and will be discussed in detail
in Section~\ref{sec:wave_equation}.

Following \cite{monk1992analysis}, we can eliminate $\bB$ from
\eqref{model:original:1}–\eqref{model:original:2} and obtain a second-order equation solely in
terms of the electric field. Taking the time derivative of \eqref{model:original:2} and using
\eqref{model:original:1} yields
\begin{equation}
\label{hyperbolc}
    \frac{\partial^2 \bE}{\partial t^2}
    + c^2(\nabla \times \nabla \times \bE)
    = -\frac{1}{\epsilon_0} \frac{d \J}{dt}.
\end{equation}
Using the vector identity
$\nabla \times \nabla \times \bE = \nabla(\nabla \cdot \bE) - \Delta \bE$ and substituting
\eqref{model:original:3} into \eqref{hyperbolc} we obtain the vector wave equation
\begin{equation}
\label{hyperbolc_2}
     \frac{\partial^2 \bE}{\partial t^2} - c^2 \Delta \bE
     = -c^2 \nabla \rho - \frac{1}{\epsilon_0} \frac{d \J}{dt}.
\end{equation}
In this work we assume that the material parameters are constant scalars (so that $c$ and
$\epsilon_0$ do not depend on space or time). Under this assumption, the Maxwell system
\eqref{model:original:1}–\eqref{model:original:4} reduces to the vector wave equation
\eqref{hyperbolc_2} for the electric field, which can be treated componentwise. This form is
particularly convenient for the space--time spectral collocation and staggered-grid discretization
developed in the remainder of the paper, since it allows us to compute $\bE$ explicitly by solving
a wave equation in space--time. The magnetic field $\bB$ is then recovered from $\bE$ in a
post-processing step using a discrete form of Faraday's law on a staggered grid (see
Sections~\ref{sec:postprocessing} and~\ref{sec:postprocessing_magnetic}). We emphasize that if
the coefficients are not constant, the decoupling leading to \eqref{hyperbolc_2} is not
straightforward; in that case one would need to work directly with the curl–curl formulation
\eqref{hyperbolc} and appropriate post-processing on staggered spaces.

Writing \eqref{hyperbolc_2} componentwise, we obtain
\begin{align}
     \frac{\partial^2 E_x}{\partial t^2}- c^2 \Delta E_x
       &= -c^2 \frac{\partial \rho}{\partial x}  - \frac{1}{\epsilon_0} \frac{d J_x}{dt},
          \label{hyperbolc_compon_1} \\
     \frac{\partial^2 E_y}{\partial t^2}- c^2 \Delta E_y
       &= -c^2 \frac{\partial \rho}{\partial y} - \frac{1}{\epsilon_0} \frac{d J_y}{dt},
          \label{hyperbolc_compon_2}\\
     \frac{\partial^2 E_z}{\partial t^2}- c^2 \Delta E_z
       &= -c^2 \frac{\partial \rho}{\partial z} - \frac{1}{\epsilon_0} \frac{d J_z}{dt},
          \label{hyperbolc_compon_3}
\end{align}
which are the scalar wave equations we discretize in space--time in
Sections~\ref{sec:collocation} and \ref{staggered:grid}.

For later use, we also rewrite \eqref{model:original:1} in decoupled form,
\begin{equation}
\label{B:decoupled}
    \begin{split}
    \frac{\partial B_x}{\partial t}
      &= -\frac{\partial E_z}{\partial y}
         + \frac{\partial E_y}{\partial z}, \\
    \frac{\partial B_y}{\partial t}
      &= +\frac{\partial E_z}{\partial x}
         - \frac{\partial E_x}{\partial z}, \\
    \frac{\partial B_z}{\partial t}
      &= -\frac{\partial E_y}{\partial x}
         + \frac{\partial E_x}{\partial y},
    \end{split}
\end{equation}
which will be the basis for reconstructing $\bB$ from the discrete electric field via the
post-processing procedures in Sections~\ref{sec:postprocessing} and~\ref{sec:postprocessing_magnetic}.

Further, we will employ the following notation throughout the text. Let
$\Omega_X,\Omega_Y,\Omega_Z \subset \mathbb{R}$ be intervals, and define the spatial computational
domain
\[
  \Omega_{\text{space}} := \Omega_X \times \Omega_Y \times \Omega_Z,
\]
so that $\Omega = \Omega_{\text{space}}$ is a cubic domain with boundary
$\partial\Omega$. The time interval is denoted by $I_T := [0,T]$, and the space--time domain is
$\Omega_T := I_T \times \Omega_{\text{space}}$. We assume such a tensor-product structure for the
computational domain in order to obtain a corresponding tensor structure for the matrices
associated with the discrete operators, which will be advantageous for tensor-train (TT)
approximations. A more detailed explanation of these tensor-network representations will be
provided in Section~\ref{sec:Tensor-Networks}.

In the remainder of this section we first introduce the space--time Chebyshev spectral collocation
scheme (Section~\ref{sec:collocation}), then construct the staggered discrete spaces for $\bE$ and
$\bB$ (Section~\ref{staggered:grid}), and finally describe the post-processing for recovering the
magnetic field from the electric field (Section~\ref{sec:postprocessing}).
\subsection{Space--time Chebyshev spectral collocation}
\label{sec:collocation}
We employ Chebyshev spectral collocation to discretize both spatial and temporal dimensions in
the Maxwell equations. This space--time formulation allows us to solve for the electric field across
all time steps simultaneously, enabling the tensor network compression techniques described in
later sections.

We apply Chebyshev polynomials $T_k(x) = \cos(k \arccos(x))$ as global basis functions.
The Chebyshev polynomials are particularly well-suited because they provide exponential
convergence for smooth solutions while maintaining favorable conditioning properties. For
computational convenience, we follow \cite{funaro1997spectral} and work with modified Chebyshev
polynomials in each coordinate direction, including time.

For each coordinate $\xi \in \{t,x,y,z\}$ we denote by $\{\xi_i\}_{i=0}^{N_\xi}$ the
Chebyshev--Gauss--Lobatto (CGL) collocation points \cite{hussaini1989spectral} on $[-1,1]$, and by
$\{l_j^{(\xi)}\}_{j=0}^{N_\xi}$ the associated interpolating polynomials, which satisfy
the property
\begin{equation*}
\label{modified}
l_j^{(\xi)}(\xi_i) =
\begin{cases}
1 & \text{if } i = j, \\
0 & \text{if } i \neq j,
\end{cases}
\end{equation*}
for every coordinate $\xi \in \{t,x,y,z\}$. To lighten notation, we will often write $l_j(\xi)$
when the underlying coordinate is clear from the context.

This construction yields interpolating polynomials based on Chebyshev nodes, which provide
the numerical advantages of Chebyshev polynomials while simplifying the implementation of
boundary conditions and the evaluation of differential operators.

The one-dimensional CGL quadrature in a generic coordinate $\xi$
provides exact integration for polynomials of degree up to $2N_\xi - 1$:
\begin{equation}
    \int_{-1}^{+1} p(\xi)\, w(\xi)\, d\xi
    = \sum_{a=0}^{N_\xi} p(\xi_a)\, \omega_a,
    \qquad
    w(\xi) = \frac{1}{\sqrt{1-\xi^2}},
\end{equation}
where $\{\omega_a\}_{a=0}^{N_\xi}$ are the Chebyshev quadrature weights.
After an affine mapping of each interval $\Omega_X,\Omega_Y,\Omega_Z, I_T$ to $[-1,1]$, we use
the same weight function $w$ in each coordinate and define the weighted continuous
$L^2$-norm on the space--time domain $\Omega_T = I_T \times \Omega_{\text{space}}$ by
\begin{equation}
    \|F\|
    := \Bigg( \int_{\Omega_T}
        |F(t,\bx)|^2\, w(t)\, w(x)\, w(y)\, w(z)\,
        \mathrm{d}t\,\mathrm{d}\bx
       \Bigg)^{1/2},
\end{equation}
where $\bx = (x,y,z)$.

Let $\{t_{a_4}\}_{a_4=0}^{N_t}$, $\{x_{a_1}\}_{a_1=0}^{N_x}$, $\{y_{a_2}\}_{a_2=0}^{N_y}$ and
$\{z_{a_3}\}_{a_3=0}^{N_z}$ denote the CGL nodes in time and the three spatial directions, with
corresponding quadrature weights $\{\omega^{(t)}_{a_4}\}$, $\{\omega^{(x)}_{a_1}\}$,
$\{\omega^{(y)}_{a_2}\}$ and $\{\omega^{(z)}_{a_3}\}$. The associated discrete norm is
\begin{equation}
    \|F\|_{N_t,N_x,N_y,N_z}
    := \Bigg(
        \sum_{a_1=0}^{N_x}
        \sum_{a_2=0}^{N_y}
        \sum_{a_3=0}^{N_z}
        \sum_{a_4=0}^{N_t}
        |F(t_{a_4},x_{a_1},y_{a_2},z_{a_3})|^2
        \,\omega^{(x)}_{a_1}\,\omega^{(y)}_{a_2}\,
        \omega^{(z)}_{a_3}\,\omega^{(t)}_{a_4}
      \Bigg)^{1/2}.
\end{equation}
In what follows, we take the same polynomial degree in each coordinate,
$N_t = N_x = N_y = N_z =: N$, and reuse the notation
$\{\xi_a\}_{a=0}^N$ and $\{\omega_a\}_{a=0}^N$ for all coordinates. In that case we simply write
\begin{equation}
\label{eq:quadformula}
    \|F\|_{N}
    := \Bigg(
        \sum_{a_1,a_2,a_3,a_4=0}^{N}
        |F(t_{a_4},x_{a_1},y_{a_2},z_{a_3})|^2
        \,\omega_{a_1}\,\omega_{a_2}\,\omega_{a_3}\,\omega_{a_4}
      \Bigg)^{1/2}.
\end{equation}
For tensor-product polynomials of degree at most $N$ in each coordinate, the discrete norm
$\|\cdot\|_{N}$ is equivalent to the continuous norm
$\|\cdot\|$. This ensures that the discrete and continuous problems have quantitatively comparable dynamics.

For later use, we also introduce the diagonal quadrature weight matrices on the one-dimensional
collocation grids. Let $W_h \in \mathbb{R}^{(N+1)\times(N+1)}$ be the diagonal matrix whose
entries are the Chebyshev weights $\{\omega_a\}_{a=0}^N$ from \eqref{eq:quadformula}.
Given a matrix $M \in \mathbb{R}^{(N+1)\times(N+1)}$ associated with a one-dimensional
collocation grid, we denote by $\langle M \rangle$ the submatrix obtained by removing the first
and last rows and columns, that is, the restriction of $M$ to the $N-1$ interior collocation points.
Similarly, for matrices acting on the three-dimensional spatial grid we write
$\langle\!\langle M \rangle\!\rangle$ for the restriction to interior spatial nodes. (The same
notation will be used later for derivative and Laplacian matrices in Section~\ref{sec:wave_equation}.)

With this convention, the tensor-product weight matrix on the interior of the space--time grid is
\begin{equation}
\label{Weight_Chebyshev}
    W := \langle\!\langle W_h \rangle\!\rangle
         \otimes \langle\!\langle W_h \rangle\!\rangle
         \otimes \langle\!\langle W_h \rangle\!\rangle
         \otimes \langle W_h \rangle.
\end{equation}
The matrix $W$ defines the weighted $\ell^2$-norm $\|W^{1/2}\theta_h\|_2$ on vectors of
interior space--time degrees of freedom, which will be used in the condition-number and
convergence analysis of Section~\ref{sec:conditioning_convergence}.

With the collocation grids, quadrature rules and norms in place, we now specify how the
Maxwell unknowns are represented in these bases. Each component of the electric and
magnetic fields is approximated by its tensor-product interpolant in the corresponding
Chebyshev--Gauss--Lobatto nodes, expanded in the interpolating polynomials introduced above.

Expanding in terms of $l_j$'s, the discrete electric and magnetic fields can be written as
\begin{equation}
\begin{split}
\label{Int_legen}
   \bE_h(t,\bx) &= (E_{xh},E_{yh},E_{zh})
   := \Bigg(
      \sum_{j=0}^{N_{E_x}} E_x^j\, l_j^{E_x}(t,\bx),\;
      \sum_{j=0}^{N_{E_y}} E_y^j\, l_j^{E_y}(t,\bx),\;
      \sum_{j=0}^{N_{E_z}} E_z^j\, l_j^{E_z}(t,\bx)
   \Bigg), \\
   \bB_h(t,\bx) &= (B_{xh},B_{yh},B_{zh})
   := \Bigg(
      \sum_{j=0}^{N_{B_x}} B_x^j\, l_j^{B_x}(t,\bx),\;
      \sum_{j=0}^{N_{B_y}} B_y^j\, l_j^{B_y}(t,\bx),\;
      \sum_{j=0}^{N_{B_z}} B_z^j\, l_j^{B_z}(t,\bx)
   \Bigg),
\end{split}
\end{equation}
where $\bx = (x,y,z)$ and the superscript on $l_j$ indicates the field component with which the
basis functions and associated collocation grid are aligned. The integers $N_{E_i}+1$ and
$N_{B_i}+1$ denote the number of spatial degrees of freedom used for each component $E_i$ and
$B_i$, respectively. These numbers differ between components due to the staggered grid
construction described in Section~\ref{staggered:grid}. This ensures proper coupling between
electric and magnetic field components through the discrete curl operator. We now construct discrete differential operators acting on these expansions.

The core of the spectral collocation method lies in constructing matrix representations of
differential operators. We build all spatial and temporal operators from the fundamental
one-dimensional derivative matrices. For the $x$-direction we define
\begin{equation}
  \left(\bold{S}_x\right)_{ij}
  := \frac{d}{dx} l_j^{(x)}(x)\Big|_{x = x_i},
  \label{time_derivative}
\end{equation}
and analogously obtain $\bS_y$, $\bS_z$ and $\bS_t$ in the $y$-, $z$-, and $t$-directions.
We construct the first-order derivative matrix $\bold{S}_x$ of size
$(N_x + 1) \times (N_x + 1)$ (and similarly for the other directions). Second-order derivative
matrices are obtained through matrix multiplication of the first-order derivatives:
\begin{equation}
\label{double:D}
  (\bold{S}_{xx})_{ij}
  = \sum_{s=0}^{N_x} (\bold{S}_x)_{is}\, (\bold{S}_x)_{sj}.
\end{equation}
A complete derivation for this expression can be found in \cite{funaro1997spectral}.

To construct operators acting on the full four-dimensional space--time domain
$\Omega_T = I_T \times \Omega_{\text{space}}$, we use Kronecker products of the
one-dimensional operators. The Laplacian operator becomes
\begin{align}
    \label{Spectral_Diff}
    \bA_{\Delta}
    &= \bold{I}_{t}\otimes \bold{S}_{xx} \otimes \bold{I}_{y}\otimes  \bold{I}_{z}
     + \bold{I}_{t} \otimes \bold{I}_{x} \otimes \bold{S}_{yy}  \otimes\bold{I}_{z}
     + \bold{I}_{t} \otimes \bold{I}_{x} \otimes  \bold{I}_{y} \otimes\bold{S}_{zz},
\end{align}
where $\bold{I}_\xi$ denotes the identity matrix of appropriate size for coordinate
$\xi \in \{x,y,z,t\}$.

The curl operator components are constructed using spatial derivatives only:
\begin{equation}
\nabla \times \mathbf{v}
= \left(
\frac{\partial v_z}{\partial y} - \frac{\partial v_y}{\partial z},\;
\frac{\partial v_x}{\partial z} - \frac{\partial v_z}{\partial x},\;
\frac{\partial v_y}{\partial x} - \frac{\partial v_x}{\partial y}
\right),
\end{equation}
where the associated matrix operators are
\begin{equation}
\label{dis:mat:1D}
\begin{split}
       & \bD_x := \bold{I}_{t}\otimes \bold{S}_{x} \otimes \bold{I}_{y}\otimes  \bold{I}_{z}, \qquad
         \bD_y := \bold{I}_{t}\otimes \bold{I}_{x} \otimes \bold{S}_{y}\otimes  \bold{I}_{z}, \\
       & \bD_z := \bold{I}_{t}\otimes \bold{I}_{x} \otimes \bold{I}_{y}\otimes  \bold{S}_{z}, \qquad
         \bD_t := \bS_{t}\otimes \bold{I}_{x} \otimes \bold{I}_{y}\otimes  \bold{I}_{z}.
\end{split}
\end{equation}
The matrices corresponding to the curl operator components are constructed as linear combinations
of $\bD_x$, $\bD_y$, and $\bD_z$, while $\bD_t$ represents the temporal derivative operator and
will be used in the reconstruction of the magnetic field.
\subsection{Staggered space--time discrete spaces}
\label{staggered:grid}
To preserve the divergence-free property of the magnetic field, $\nabla \cdot \mathbf{B} = 0$, in our discretization, we introduce a staggered spectral collocation approach within the space--time framework developed in Section~2.1.

In classical finite-difference staggered grids, different field components are evaluated at spatially shifted locations to improve stability and enforce discrete conservation laws. In our spectral method, we achieve analogous staggering not through spatial shifts, but by using \emph{different polynomial degrees} for electric and magnetic field components in different coordinate directions.

Specifically, we collocate the electric field components $\mathbf{E}$ at Chebyshev--Gauss--Lobatto (CGL) nodes, which correspond to polynomials of a certain degree, while magnetic field components $\mathbf{B}$ are collocated at Chebyshev--Gauss (CG) nodes---the interior points that exclude the boundaries. This spectral staggering maintains exponential convergence while inheriting the stability and conservation properties of classical staggered grids.

The key insight is that the discrete curl operator $\nabla_h \times \mathbf{E}_h$, assembled from the one-dimensional derivative matrices in~\eqref{dis:mat:1D}, must map between compatible discrete spaces. By carefully constructing polynomial spaces with different degrees of freedom for each field component, we ensure that:
\begin{enumerate}
\item The discrete curl accurately approximates $\partial \mathbf{B}_h/\partial t$ at each time step
\item The discrete identity $\nabla \cdot (\nabla \times \cdot) = 0$ holds to spectral accuracy
\item The constraint $\nabla_h \cdot \mathbf{B}_h \approx 0$ is maintained throughout the computation
\end{enumerate}

Unlike finite-difference methods where staggering involves pointwise evaluations at shifted spatial locations, our spectral staggered grid uses different polynomial degrees for electric and magnetic field components in different coordinate directions~\cite{niu2023legendre,huang2012divergence}. This approach preserves spectral accuracy while providing the structural benefits of classical staggering.

Recall that the one-dimensional computational domains in each spatial
direction are intervals $\Omega_X,\Omega_Y,\Omega_Z \subset \mathbb{R}$ and
that the time interval is $I_T = [0,T]$, so that
$\Omega_{\text{space}} = \Omega_X \times \Omega_Y \times \Omega_Z$ and
$\Omega_T = I_T \times \Omega_{\text{space}}$.
Under the constant-coefficient assumptions stated in Section~\ref{sec:math:model:num:dis}, the Maxwell system reduces to the vector wave equation~\eqref{hyperbolc_2} for the electric field.
The design of the discrete spaces is crucial to ensure that
the discrete curl and divergence operators satisfy the same structural
relations as in the continuous case.

The key insight is that the discrete curl operator
$\nabla_h \times \bE_h$, assembled from the one-dimensional derivative
matrices in~\eqref{dis:mat:1D}, must accurately approximate
$\partial \bB_h / \partial t$ at each time step and map between compatible
discrete spaces. This requires a careful construction of polynomial spaces
with different degrees of freedom for each field component, so that the
discrete analogue of the identity $\nabla \cdot (\nabla \times \cdot) = 0$
holds asymptotically and the constraint $\nabla_h \cdot \bB_h \approx 0$
is maintained.

Let $\mathbb{P}_k(I)$ denote the space of polynomials of degree at most $k$
on an interval $I$. We construct distinct tensor-product polynomial spaces for
each component of the magnetic field
$\mathbf{B} = (B_x, B_y, B_z)$:
\begin{equation}
\label{grid:B}
    \begin{split}
        S^x_k &:= \mathbb{P}_{k+1}(I_T) \otimes \mathbb{P}_{k+1}(\Omega_X)
                 \otimes \mathbb{P}_k(\Omega_Y) \otimes \mathbb{P}_k(\Omega_Z),\\
        S^y_k &:= \mathbb{P}_{k+1}(I_T) \otimes \mathbb{P}_{k}(\Omega_X)
                 \otimes \mathbb{P}_{k+1}(\Omega_Y) \otimes \mathbb{P}_k(\Omega_Z), \\
        S^z_k &:= \mathbb{P}_{k+1}(I_T) \otimes \mathbb{P}_{k}(\Omega_X)
                 \otimes \mathbb{P}_{k}(\Omega_Y) \otimes \mathbb{P}_{k+1}(\Omega_Z).
    \end{split}
\end{equation}
For the magnetic field component $B_i$, we use polynomials of degree $k + 1$ in the temporal direction and in the corresponding $i$-th spatial direction, and polynomials of degree $k$ in the remaining two spatial directions. In practice, this means $B_x$ is collocated at $(N+1) \times (N+1) \times N \times N$ points in the $(t,x,y,z)$ directions, while $E_x$ is collocated at $(N+1) \times N \times (N+1) \times (N+1)$ points, where $N=k+1$.

For the electric field components
$\mathbf{E} = (E_x, E_y, E_z)$, we define complementary spaces:
\begin{equation}
\label{grid:E}
    \begin{split}
        V^x_k &:= \mathbb{P}_{k+1}(I_T) \otimes \mathbb{P}_{k}(\Omega_X)
                 \otimes \mathbb{P}_{k+1}(\Omega_Y) \otimes \mathbb{P}_{k+1}(\Omega_Z),\\
        V^y_k &:= \mathbb{P}_{k+1}(I_T) \otimes \mathbb{P}_{k+1}(\Omega_X)
                 \otimes \mathbb{P}_{k}(\Omega_Y) \otimes \mathbb{P}_{k+1}(\Omega_Z), \\
        V^z_k &:= \mathbb{P}_{k+1}(I_T) \otimes \mathbb{P}_{k+1}(\Omega_X)
                 \otimes \mathbb{P}_{k+1}(\Omega_Y) \otimes \mathbb{P}_{k}(\Omega_Z).
    \end{split}
\end{equation}
Here, the electric field component $E_i$ uses degree $k$ in the $i$-th
spatial direction and degree $k+1$ in the other two spatial directions, while
being of degree $k+1$ in time. This choice is complementary to the
construction of the magnetic-field spaces in~\eqref{grid:B}.

This staggered construction ensures that when we compute
\[
\frac{\partial \bB_h}{\partial t}
= -\,\nabla_h \times \bE_h
\]
using the discrete curl operator built from
$\bD_x,\bD_y,\bD_z$ in~\eqref{dis:mat:1D}, the resulting magnetic-field
components naturally lie in the spaces $S^i_k$ and satisfy the discrete
divergence-free constraint to high accuracy. In particular, the curl
operation maps from the electric field spaces $V^i_k$ to the magnetic field
spaces $S^i_k$ while preserving spectral accuracy and maintaining the proper
relationships between degrees of freedom.

The choice of which components use degree $k$ versus $k+1$ in each direction is dictated by the structure of Maxwell's curl equations—ensuring that when we differentiate a degree $k+1$ polynomial, we obtain a degree $k$ polynomial that naturally lives in the correct discrete space.
We now collect the component-wise spaces into vector-valued discrete spaces
for the magnetic and electric fields:
\begin{equation}
\label{Vspaces}
\begin{split}
\bScal_h &:= S^x_k \times S^y_k \times S^z_k
\quad \text{(for the magnetic field $\bB_h$)}, \\
\bV_h &:= V^x_k \times V^y_k \times V^z_k
\quad \text{(for the electric field $\bE_h$)}.
\end{split}
\end{equation}
Thus, each component $B_i$ (respectively, $E_i$) belongs to $S^i_k$
(respectively, $V^i_k$), and the expansions in~\eqref{Int_legen} are
consistent with these spaces and the associated staggered collocation grids.

At the level of the wave equation, the fully discrete space--time scheme
for the electric field can be written formally as: find
$\bE_h \in \bV_h$ such that
\begin{equation}
\frac{\partial^2 \bE_h}{\partial t^2}
  - c^2 \Delta \bE_h
  = -c^2 \nabla \rho_h
    - \frac{1}{\epsilon_0} \frac{d\J_h}{dt},
\label{hyperbolc_discrete}
\end{equation}
where $\rho_h$ and $\J_h$ denote the interpolants of the charge and current
densities onto the collocation grid, and the differential operators are
discretized using the matrices introduced in Section~\ref{sec:collocation}
(e.g.\ the Laplacian matrix $\bA_\Delta$ in~\eqref{Spectral_Diff} and the
first-order derivative matrices in~\eqref{dis:mat:1D}). The explicit matrix
form of~\eqref{hyperbolc_discrete}, including the Kronecker structure of the
resulting operator, will be derived in Section~3.

The different degrees of freedom for each component in the staggered spaces
$\bS_h$ and $\bV_h$ require a dedicated post-processing step to recover the
magnetic field from the electric field solution. This post-processing,
based on the discrete Faraday law, is outlined in
Section~\ref{sec:postprocessing} and later revisited in matrix form in
Section~3.1.
\subsection{Post-processing to recover the magnetic field}
\label{sec:postprocessing}

Our space--time spectral discretization solves the wave equation for the electric field and
computes a space--time approximation $\bE_h \in \bV_h$ across all collocation points
(see Sections~\ref{sec:collocation} and~\ref{staggered:grid}). The magnetic field
$\bB_h \in \bScal_h$ is then obtained in a separate, but fully consistent, post-processing step
based on the discrete version of Faraday’s law:
\begin{equation}
\label{eq:faraday_discrete}
    \frac{\partial \bB_h}{\partial t} = - \nabla_h \times \bE_h.
\end{equation}
Here $\nabla_h \times$ denotes the discrete curl operator assembled from the one-dimensional
derivative matrices introduced in~\eqref{dis:mat:1D}. This equation is evaluated entirely within
the tensor-product spectral framework described in
Sections~\ref{sec:collocation}--\ref{staggered:grid}.

We view the discrete electric field $\bE_h$ given by the space--time solve in its interpolant
form based on the CGL nodes. Each component of the electric field can be
written as
\begin{equation}
\label{eq:Eh_interp}
    E_{i h }(t,\bx)
    = \sum_{a_1,a_2,a_3,a_4=0}^{N}
      E_{i}^{a_1,a_2,a_3,a_4} \;
      l^{E_i}_{a_1}(x)\, l^{E_i}_{a_2}(y)\, l^{E_i}_{a_3}(z)\, l^{E_i}_{a_4}(t),
    \qquad i \in \{x,y,z\},
\end{equation}
where $\bx = (x,y,z)$, the functions $l^{E_i}_{a_j}$ are the component-wise basis
functions associated with the staggered collocation grids, and
$E_{i}^{a_1,a_2,a_3,a_4}$ are the coefficients obtained from the space--time solve in
$\bV_h$. This expression particularly is very useful to derive the low-rank structure of the associated matrices.

To compute the magnetic field from~\eqref{eq:faraday_discrete}, we apply the discrete curl
operator built from the matrices in~\eqref{dis:mat:1D} to the interpolant~\eqref{eq:Eh_interp}.
By construction of the staggered spaces in Section~\ref{staggered:grid}, the electric field
lives in $\bV_h$ and the discrete curl $\nabla_h \times$ maps from $\bV_h$ into
$\bScal_h$. Consequently, the reconstructed magnetic field $\bB_h$ automatically belongs to the
correct discrete space and remains compatible with the electric field representation.
This staggered mapping preserves the discrete divergence-free condition
$\nabla_h \cdot \bB_h \approx 0$ and maintains the overall spectral accuracy of the scheme.

The discussion in this section is intentionally kept at a conceptual level. In
Section~3.1, we revisit the magnetic field recovery and derive the corresponding
matrix equations in detail, highlighting the Kronecker and tensor-train structure of
the post-processing step.
\section{Wave equation in the space--time spectral collocation scheme}
\label{sec:wave_equation}

In Sections~\ref{sec:collocation}--\ref{sec:postprocessing}, we introduced the space--time
Chebyshev spectral collocation framework and the staggered spaces
$\bV_h$ and $\bScal_h$ for the electric and magnetic fields, respectively. Under the constant-coefficient
assumptions in Section~\ref{sec:math:model:num:dis}, the Maxwell system reduces to the vector wave
equation~\eqref{hyperbolc_2} for the electric field. In this section, we derive the explicit space--time
matrix formulation of this wave equation for a single component of the electric field. This will reveal
the Kronecker structure of the discrete operator and motivate the low-rank tensor representations
exploited in later sections.

For notational convenience we abbreviate the spatial domain by
$\Omega := \Omega_{\text{space}}$, so that $\Omega_T = I_T \times \Omega$.
We detail the derivation for the $x$-component of the electric field; the $y$- and $z$-components
are treated analogously and lead to systems of the same form.

We begin by rewriting the wave equation \eqref{hyperbolc_compon_1} in first-order form and then
present its discretization using the space--time spectral collocation method introduced in
Section~\ref{sec:collocation}. The resulting system naturally admits a low-rank tensor
representation, which will be exploited in later sections.

Let us define the auxiliary variable for the $x$-component of the electric field as
\[
\mathbf{v}
= \begin{bmatrix} v \\ v_t \end{bmatrix}
:= \begin{bmatrix} E_x \\ \frac{\partial E_x}{\partial t} \end{bmatrix},
\]
and define the corresponding source term as
\[
f := -c^2 \,\frac{\partial \rho}{\partial x} - \frac{1}{\epsilon_0} \frac{d J_x}{dt}.
\]
Then, the second-order wave equation can be written as the first-order system
\begin{equation}
\frac{d\mathbf{v}}{dt}
=
\begin{bmatrix}
0 & I \\
c^2 \Delta & 0
\end{bmatrix} \mathbf{v} +
\begin{bmatrix}
0 \\
f
\end{bmatrix},
\label{eq:wave_first_order}
\end{equation}
with boundary and initial conditions
\[
\mathbf{v}(t, \partial \Omega) =
\begin{bmatrix}
v^{\text{bd}}(t, \mathbf{x}) \\
v^{\text{bd}}_t(t, \mathbf{x})
\end{bmatrix},
\qquad
\mathbf{v}(0, \Omega) =
\begin{bmatrix}
v^0(\mathbf{x}) \\
v_t^0(\mathbf{x})
\end{bmatrix},
\]
where $v^{\text{bd}}$ and $v^{\text{bd}}_t$ prescribe boundary data for $E_x$ and its time derivative on
$\partial\Omega$, while $v^0$ and $v_t^0$ denote the initial value and the initial time
derivative of $E_x$ at $t=0$.

Let $\mathbf{v}_h$ denote the values of $\mathbf{v}$ evaluated at the space--time collocation
grid. The discretized form of \eqref{eq:wave_first_order} becomes
\begin{equation}
    \begin{bmatrix}
        \bS_t \otimes \bI_{(N+1)^3}  & \bZero \\
        \bZero & \bS_t \otimes \bI_{(N+1)^3}
    \end{bmatrix} 
    \begin{bmatrix}
        \bv_{h} \\
        \bv_{th}
    \end{bmatrix}
    =
    \begin{bmatrix}
        \bZero  & \bI_{N+1}\otimes \bI_{(N+1)^3} \\
        \bI_{N+1} \otimes \bA_{\Delta} & \bZero
    \end{bmatrix} 
    \begin{bmatrix}
        \bv_h \\
        \bv_{th}
    \end{bmatrix}
    +
    \begin{bmatrix}
        \bZero \\
        \bff_h
    \end{bmatrix},
    \label{wave:1}
\end{equation}
where $\mathbf{f}_h$ is the discretized source, and $\bI_m$ denotes the $m \times m$ identity
matrix. The matrix $\bA_\Delta$ is the discrete Laplacian defined in~\eqref{Spectral_Diff}. 

We next eliminate the known boundary and initial conditions and restrict the system to interior
space--time points. Let $\hat{\mathbf{v}}_{h}$ and $\hat{\mathbf{v}}_{th}$ denote the unknown
interior values of $E_x$ and its time derivative. Denoting by
$\langle\bS_t\rangle$ and $\langle\!\langle \bA_{\Delta} \rangle\!\rangle$ the restrictions of
$\bS_t$ and $\bA_\Delta$ to interior time and space nodes, respectively, we obtain
\begin{equation}
\begin{split}
    \begin{bmatrix}
        \langle\bS_t \rangle \otimes \bI_{(N-1)^3}  & \bZero \\
        \bZero & \langle \bS_t \rangle \otimes \bI_{(N-1)^3}
    \end{bmatrix} 
    \begin{bmatrix}
        \hat{\bv}_{h} \\
        \hat{\bv}_{th}
    \end{bmatrix}
    &=
    \begin{bmatrix}
        \bZero  & \bI_{N}\otimes \bI_{(N-1)^3} \\
       \bI_{N} \otimes \langle\!\langle \bA_{\Delta} \rangle\!\rangle & \bZero
    \end{bmatrix} 
    \begin{bmatrix}
        \hat{\bv}_{h} \\
        \hat{\bv}_{th}
    \end{bmatrix}
   \\
  & \quad +\begin{bmatrix}
        \bZero \\
        \hat{\bff}_h
    \end{bmatrix}
    -\begin{bmatrix}
        \bA_t \bv^{0}_{h} \\
        \bA_t \bv^{0}_{th}
    \end{bmatrix}
    + \begin{bmatrix}
        \mathbf{G} \bv^{\text{bd}}_{h} \\
        \mathbf{G}_t \bv^{\text{bd}}_{th}
    \end{bmatrix},
     \label{wave:2}
\end{split}
\end{equation}
where
\begin{equation}
\begin{split}
    \bA_t &:= \bS_t(I^{\text{int}}_t,:) \otimes \bold{I}_z(I^{\text{int}}_z,:) \otimes \bold{I}_y(I^{\text{int}}_y,:) \otimes \bold{I}_x(I^{\text{int}}_x,:), \\
    \mathbf{G} &:= \bold{I}_t(I^{\text{int}}_t,:) \otimes \bold{I}_z(I^{\text{int}}_z,:) \otimes \bold{I}_y(I^{\text{int}}_y,:) \otimes \bold{I}_x(I^{\text{int}}_x,:), \\
    \mathbf{G}_t &:= \bold{I}_t(I^{\text{int}}_t,:) \otimes \bold{I}_z(I^{\text{int}}_z,:) \otimes \bold{I}_y(I^{\text{int}}_y,:) \otimes \bS_{xx}(I^{\text{int}}_x,:) \\
          &\qquad + \bold{I}_t(I^{\text{int}}_t,:) \otimes \bold{I}_z(I^{\text{int}}_z,:) \otimes \bS_{yy}(I^{\text{int}}_y,:) \otimes \bold{I}_x(I^{\text{int}}_x,:) \\
          & \qquad + \bold{I}_t(I^{\text{int}}_t,:) \otimes \bS_{zz}(I^{\text{int}}_z,:) \otimes \bold{I}_y(I^{\text{int}}_y,:) \otimes \bold{I}_x(I^{\text{int}}_x,:).
\end{split}
\end{equation}
Here, for a matrix $\bold{M}$ and an index set $I$, the notation $\bold{M}(I,:)$ denotes the submatrix
consisting of rows indexed by $I$. The sets $I^{\text{int}}_t, I^{\text{int}}_x, I^{\text{int}}_y, I^{\text{int}}_z$ collect the indices of interior nodes in time and in each spatial direction, respectively. The boundary data
is stored in the vectors $\bv^{\text{bd}}_{h}$ and $\bv^{\text{bd}}_{th}$ (discrete counterparts
of $v^{\text{bd}}$ and $v^{\text{bd}}_t$), initial data in vectors $\bv^0_{h}$ and $\bv^0_{th}$ (discrete counterparts
of $v^0$ and $v^0_t$), and
$\hat{\bff}_h$ denotes the interior part of the source after boundary elimination.

Decoupling \eqref{wave:2} we get:
\begin{equation}
    \hat{\bv}_{th}
    =(\langle \bS_t \rangle \otimes \bI_{(N-1)^3})\, \hat{\bv}_{h}
      + \bA_t \bv^0_{h}
      - \mathbf{G} \bv^{\text{bd}}_{h},
    \label{v:1}
\end{equation}
from the first block and the following in the second block
\begin{equation}
    (\bI_N \otimes \langle\!\langle \bA_{\Delta} \rangle\!\rangle)\hat{\bv}_{h}
    = (\langle \bS_t \rangle \otimes \bI_{(N-1)^3})\hat{\bv}_{th}
      - \hat{\bff}_h
      +\bA_t \bv^0_{th}
      -\mathbf{G}_t \bv^{\text{bd}}_{th}.
    \label{v:2}
\end{equation}

Eliminating $\hat{\bv}_{th}$, we get a system in terms of $\hat{\bv}_h$:
\begin{equation}
\bA_{Lap}\, \hat{\mathbf{v}}_{h} =  \mathbf{F} - \mathbf{F}^{\text{BD}},
\label{eq:wave_final}
\end{equation}
with
\begin{align*}
    \bA_{Lap} &:= \big( \langle \bS_t \rangle^2 \otimes \bI_{(N-1)^3} \big)
                 - \big( \bI_N \otimes \langle\!\langle \bA_{\Delta} \rangle\!\rangle \big),\\
    \mathbf{F} &= \hat{\bff}_h, \\
    \mathbf{F}^{\text{BD}} &=
        \bA_t \bv^0_{th}
        -\mathbf{G}_t \bv^{\text{bd}}_{th}
        +(\langle \bS_t \rangle \otimes \bI_{(N-1)^3}) \bA_t \bv^0_{h}
        - (\langle \bS_t \rangle \otimes \bI_{(N-1)^3}) \mathbf{G} \bv^{\text{bd}}_{h},
\end{align*}
where $\mathbf{F}^{\text{BD}}$ collects all the discrete initial and boundary contributions. Thus, the fully discrete space--time wave equation
for the $x$-component reduces to the linear system~\eqref{eq:wave_final}. The same discretization
procedure applies to the $y$- and $z$-components of the electric field. Replacing $E_x$ by $E_y$ or
$E_z$, and $J_x$ by $J_y$ or $J_z$, respectively, leads to systems of the form \eqref{eq:wave_final}
for each spatial direction. The tensor structure and solution strategy remain identical.

The operator $\bA_{Lap}$ in \eqref{eq:wave_final} is assembled using Kronecker products of
one-dimensional matrices in time and each spatial dimension. This structure is ideal for
approximation in TT format, enabling scalable algorithms in higher dimensions.

\subsection{Recovering the magnetic field}
\label{sec:postprocessing_magnetic}

Once the electric field components
$\bE_h = (\bE_{xh}, \bE_{yh}, \bE_{zh}) \in \bV_h$ have been computed from the space--time
wave equation~\eqref{eq:wave_final}, we recover the magnetic field
$\bB_h = (\bB_{xh}, \bB_{yh}, \bB_{zh}) \in \bScal_h$ using the discrete version of Faraday's law, as
previewed in Section~\ref{sec:postprocessing}. The discrete
magnetic field components are governed by the system
\begin{equation}
\label{B:decoupled:decoupled}
    \begin{split}
    \frac{\partial B_{xh}}{\partial t}
      &= \underbrace{-\frac{\partial E_{zh}}{\partial y}
                     + \frac{\partial E_{yh}}{\partial z}}_{=:H_{xh}}, \\
    \frac{\partial B_{yh}}{\partial t}
      &= \underbrace{+\frac{\partial E_{zh}}{\partial x}
                     - \frac{\partial E_{xh}}{\partial z}}_{=:H_{yh}}, \\
    \frac{\partial B_{zh}}{\partial t}
      &= \underbrace{-\frac{\partial E_{yh}}{\partial x}
                     + \frac{\partial E_{xh}}{\partial y}}_{=:H_{zh}}.
    \end{split}
\end{equation}
Here $H_{ih}$ represent the curl components of the electric field. Using the matrix operators
defined in~\eqref{dis:mat:1D}, each magnetic field component satisfies a linear system of the form
\begin{equation}
\bD_t \bB_{ih} = \bH_{ih}, \quad i \in \{x,y,z\}, \label{eq:magnetic_system}
\end{equation}
where $\bD_t := \bS_t \otimes \bI_x \otimes \bI_y \otimes \bI_z$ is the temporal derivative
operator and $\bB_{ih}, \bH_{ih}$ denote the vectors of collocation values of $B_{ih}$ and $H_{ih}$.

For the $x$-component of the magnetic field, the right-hand side $\bH_{xh}$ is computed using the
spatial derivative operators:
\begin{equation}
\bH_{xh}
= -(\bI_t \otimes \bI_x \otimes \bS_y \otimes \bI_z)\, \bE_{zh}
  + (\bI_t \otimes \bI_x \otimes \bI_y \otimes \bS_z)\, \bE_{yh}.
\label{eq:H1_computation}
\end{equation}
Due to the staggered grid construction described in Section~\ref{staggered:grid}, the electric
field components $E_{yh}$ and $E_{zh}$ belong to the spaces $S_k^y$ and $S_k^z$, respectively,
while the magnetic field component $B_{xh}$ belongs to the space $V_k^x$. The interpolation
between these spaces is handled through projection matrices that preserve spectral accuracy.

After applying boundary conditions, the system \eqref{eq:magnetic_system} for interior nodes
becomes
\begin{equation}
\bA_{\text{curl}} \hat{\bB}_{xh} = \bF_{\text{curl}} - \bF_{\text{curl}}^{\text{BD}},
\label{eq:magnetic_interior}
\end{equation}
where
\begin{align}
\bA_{\text{curl}} &:= \langle \bS_t \rangle \otimes \bI_x^{\text{int}} \otimes \bI_y^{\text{int}} \otimes \bI_z^{\text{int}}, \\
\bF_{\text{curl}} &:= \hat{\bH}_{xh}, \label{curlC}\\
\bF^{\text{BD}}_{\text{curl}} &:= \bD_t^{\text{int}} \bB_{xh}^{\text{bd}}.
\end{align}

In these expressions, the superscript $\text{int}$ denotes restriction to interior collocation
points. For each spatial coordinate $\xi \in \{x,y,z\}$ we denote by $I_\xi^{\text{int}}$ the
index set of interior nodes in that direction, and we define the corresponding interior identity
matrices as
\[
    \bI^{\text{int}}_\xi := \bI_\xi(I_\xi^{\text{int}},\, I_\xi^{\text{int}}),
    \qquad \xi\in\{x,y,z\}.
\]
Similarly, the temporal–spatial derivative operator restricted to the interior is defined by
\[
    \bD_t^{\text{int}}
    := \bS_t(I_t^{\text{int}}, :)
       \,\otimes\, \bI_x(I_x^{\text{int}}, :)
       \,\otimes\, \bI_y(I_y^{\text{int}}, :)
       \,\otimes\, \bI_z(I_z^{\text{int}}, :).
\]
  The vector
$\boldsymbol{B}_{xh}^{\text{bd}}$ contains the boundary values of $B_{xh}$, and boldface symbols
such as $\bB_{xh}$ and $\bF_{\text{curl}}$ represent vectors of discrete
collocation values rather than continuous fields.

The curl computation in~\eqref{curlC} requires interpolation between different staggered spaces.
Specifically, we have
\begin{align}
\bF_{\text{curl}}
= &-(\bI_t^{\text{int}} \otimes \bI_x^{\text{int}} \otimes \langle \langle \bS_y \rangle \rangle \otimes \bI_z^{\text{int}})
     \,\boldsymbol{\mathcal{I}}_{V_z \to S_x}\, \hat{\bE}_{zh} \nonumber \\
  &+(\bI_t^{\text{int}} \otimes \bI_x^{\text{int}} \otimes \bI_y^{\text{int}} \otimes \langle \langle \bS_z \rangle \rangle)
     \,\boldsymbol{\mathcal{I}}_{V_y \to S_x}\, \hat{\bE}_{yh},
\label{eq:F_curl_detailed}
\end{align}
where $\boldsymbol{\mathcal{I}}_{V_z \to S_x}$ and $\boldsymbol{\mathcal{I}}_{V_y \to S_x}$
represent interpolation operators mapping from electric-field spaces $V_k^z$ and $V_k^y$ to the
magnetic-field space $S_k^x$. The magnetic field recovery for components $B_{yh}$ and $B_{zh}$
follows the same pattern as described above for $B_{xh}$, with appropriate modifications to the
spatial derivative operators and interpolation matrices according to the staggered grid
construction in \eqref{grid:B}--\eqref{grid:E}.

This post-processing approach ensures that the reconstructed magnetic field $\bB_h$ maintains the
same spectral accuracy as the electric field solution, preserves the discrete divergence-free
constraint $\nabla_h \cdot \bB_h \approx 0$, and maintains the tensor-product structure required
for efficient implementation in tensor-train format.
\section{Condition number estimates and spectral convergence}
\label{sec:conditioning_convergence}

In this section we study two complementary analytical properties of the proposed space--time
spectral collocation scheme. First, we derive condition number estimates for the discrete
operators arising in the wave equation formulation and in the magnetic-field reconstruction,
namely the matrices $\bA_{Lap}$ and $\bA_{\text{curl}}$ introduced in
Sections~\ref{sec:wave_equation} and~\ref{sec:postprocessing_magnetic}. These estimates
quantify the stability of the discrete systems and are relevant for the performance of iterative
solvers and tensor-train approximations. Second, we establish spectral convergence of the
discrete electric and magnetic fields, as well as exponential decay of the discrete divergence
errors associated with Gauss's laws.

\subsection{Condition number estimates}

We recall the spectral condition number of a matrix $\bM$, defined by
\begin{equation}
    \kappa(\bM)
    := \frac{\max_{\lambda \in \Lambda(\bM)} |\lambda|}
            {\min_{\lambda \in \Lambda(\bM)} |\lambda|},
\end{equation}
where $\Lambda(\bM)$ denotes the spectrum of $\bM$. Throughout this subsection, for a complex number $z$ we write $Re(z)$ and $Im(z)$
for the real and imaginary parts of $z$, respectively. 

We collect several lemmas that describe
the spectral properties of Chebyshev derivative matrices \cite{lui2021spectral}. Here $\bS$ denotes a generic
one-dimensional Chebyshev first-derivative matrix (in any of the variables $t,x,y,z$), and
$\langle \bS \rangle$, $\langle\!\langle \bS^2 \rangle\!\rangle$ denote the corresponding
restrictions to interior nodes, as in Section~\ref{sec:wave_equation}.

\begin{lemma}
\label{Lemma:new}
    Let $N \geq 1$, and $\lambda \in \Lambda (\langle \bS \rangle)$. Then
    \[
        Re(\lambda) \geq C N,
    \]
    where $C$ is a positive constant independent of $N$.
\end{lemma}

\begin{lemma}
\label{Lemma:1}
    Let $N \geq 1$, and $\lambda \in \Lambda (\langle \bS \rangle)$. Then
    \[
        |\lambda| \leq c N^2,
    \]
    where $c$ is a positive constant independent of $N$.
\end{lemma}

\begin{lemma}
\label{Lemma:2}
    Let $N \geq 2$. Then the eigenvalues of $- \langle\!\langle \bS^2 \rangle\!\rangle$ are real,
    bounded below by $c$ and above by $C N^4$, where $c$ and $C$ are positive constants
    independent of $N$.
\end{lemma}

\begin{theorem}
    Let $N\geq 2$, and let $\bA_{Lap}$ be the Chebyshev spectral collocation matrix defined in
    \eqref{eq:wave_final}, and $\lambda$ be any eigenvalue of $\bA_{Lap}$. Then
    \begin{equation}
        c \leq |\lambda| \leq C N^4.
    \end{equation}
    Consequently,
    \begin{equation}
        \kappa(\bA_{Lap}) \leq C N^{4}.
    \end{equation}
\end{theorem}

\begin{proof}
    Let $\lambda \in \Lambda(\bA_{Lap})$, where $\Lambda(\bA_{Lap})$ denotes the spectrum of
    $\bA_{Lap}$. Then we can write
    \begin{equation}
        \lambda = \gamma^2 + 3 \mu,
    \end{equation}
    where $\gamma = Re(\gamma) + i\,Im(\gamma)$ is an eigenvalue of $\langle \bS_t \rangle$, and
    $\mu$ is an eigenvalue of $- \langle\!\langle \bS_{\nu}^2 \rangle\!\rangle$, with
    $\nu \in \{x,y,z\}$. A straightforward calculation yields
\begin{equation}
\begin{split}
    |\lambda|^2
      &= Re(\gamma)^4 + 6 \mu\,Re(\gamma)^2
         + 2 Re(\gamma)^2 Im(\gamma)^2
         + (3\mu - Im(\gamma))^2 \\
      &\geq Re(\gamma)^4 + 6 \mu\,Re(\gamma)^2
       \;\geq\; c,
\end{split}
\end{equation}
    using Lemma~\ref{Lemma:new} and Lemma~\ref{Lemma:2}. Moreover, by Lemma~\ref{Lemma:1} we
    obtain
    \begin{equation}
        |\lambda|^2 \leq C N^8.
    \end{equation}
    Hence $c \leq |\lambda|\leq C N^4$, and therefore $\kappa(\bA_{Lap}) \leq C N^4$.
\end{proof}

We now focus on the spectrum of the matrix associated with the magnetic-field reconstruction
defined in~\eqref{eq:magnetic_interior}.

\begin{theorem}
    Let $N \geq 1$, and let $\bA_{\text{curl}}$ be the Chebyshev spectral collocation matrix defined
    in \eqref{eq:magnetic_interior}, and $\lambda$ be an eigenvalue of $\bA_{\text{curl}}$. Then
    \begin{equation}
      c \leq |\lambda| \leq C N^2,
    \end{equation}
    and
    \begin{equation}
        \kappa(\bA_{\text{curl}}) \leq C N^2,
    \end{equation}
    where $C$ is a positive generic constant independent of $N$.
\end{theorem}

\begin{proof}
    The result follows directly from Lemma~\ref{Lemma:new} and Lemma~\ref{Lemma:1}, applied
    to the temporal derivative matrix appearing in $\bA_{\text{curl}}$.
\end{proof}

\subsection{Spectral convergence}

We now discuss the spectral convergence of the fully discrete scheme for both the wave
equation and the Maxwell system. The first result concerns the componentwise wave equation
for the electric field; the second addresses the magnetic field obtained via post-processing; and
the third quantifies the convergence of the discrete Gauss law for $\bE$ and the discrete
divergence-free constraint for $\bB$.

\begin{theorem}\label{Th:E}
    Let $E_x$ be the solution of \eqref{hyperbolc_compon_1}. Assume that $E_x$ is an analytic
    function in $t,x,y$, and $z$. Let $N \geq 2$, and let $\hat{v}_{h}$ be the solution of the
    space--time method defined in Section~\ref{sec:wave_equation}. Further, let $\Theta_h$ be the
    error associated with the collocation approximation $\hat{v}_{h}$ of $E_x$. Then
    \begin{equation}
        \bigl| W^{1/2}\Theta_h \bigr|_2 \leq c\, N ^{4.5} e^{-C N},
    \end{equation}
    where the weight matrix $W$ is defined in \eqref{Weight_Chebyshev}, and $|\cdot|_2 $ denotes vector/ matrix 2 norm.
\end{theorem}

\begin{proof}
        We assume that $\Omega_h$ is the computational grid. Let $N^{\mathrm{dof}}$ be the total
        number of grid points, and $N^{\mathrm{dof}}_{\mathrm{int}}$ be the total number of interior
        points in space. Then, we can write $E_{xh}(t)$ as
        \begin{equation*}
            E_{xh}(t):=
            \begin{bmatrix}
               E_x(t,\bx_1) \\
               E_x(t,\bx_2) \\
               \vdots \\
               E_x(t,\bx_{N^{\mathrm{dof}}_{\mathrm{int}}})
            \end{bmatrix}_{N^{\mathrm{dof}}_{\mathrm{int}} \times 1},
        \end{equation*}
        and similarly
        \begin{equation*}
            f_h(t):=
            \begin{bmatrix}
               f(t,\bx_1) \\
               f(t,\bx_2) \\
               \vdots \\
               f(t,\bx_{N^{\mathrm{dof}}_{\mathrm{int}}})
            \end{bmatrix}_{N^{\mathrm{dof}}_{\mathrm{int}} \times 1}.
        \end{equation*}
        A semi-discrete approximation of the wave equation \eqref{hyperbolc_compon_1} is 
        \begin{equation}
        \begin{split}
           \frac{d^2 E_{xh}(t)}{d t^2}
           &= \sum_{j=2}^N \Big ( \bA\,E_{xh}(t_j)+ f_h(t_j) \Big ) l_j (t), \\
           &\qquad E_{xh}(0, \Omega)=v^0,\quad \frac{d E_{xh}}{dt}(0, \Omega)= v^0_t, \\
           &\qquad E_{xh}(t, \partial \Omega)=v^{\text{bd}},\quad \frac{d E_{xh}}{dt}(t,\partial \Omega)= v^{\text{bd}}_t,
        \end{split}
        \end{equation}
        where $\bA = \langle\!\langle \bA_{\Delta} \rangle\!\rangle$. Hence, we have the equality
        \begin{equation}
        \label{error:1}
            \frac{d^2 E_{xh}(t_k)}{d t^2}
            =  \bA\,E_{xh}(t_k)+ f_h(t_k),
            \qquad 2 \leq k \leq N.
        \end{equation}
        Define the error function $e_h(t)=E_{xh}(t)-E_x(t,\bx_h)$ with $j-th$ components
        $e_j(t):= (e_h(t))_j$, where $\bx_h$ is a grid point. Using \eqref{error:1} and
        \eqref{hyperbolc_compon_1}, we conclude that the error satisfies, for $2 \leq k \leq N$,
        \begin{equation}
            \frac{d^2 e_h(t_k)}{d t^2}=  \bA\, e_h(t_k)+ r(t_k),
            \qquad
            r(t_k)= \bA E_x(t_k,\bx_h)-\Delta E_{x}(t_k,\bx_h).
        \end{equation}
We define \begin{equation*}
            \Theta_h:=
            \begin{bmatrix}
               e_1(t_h) \\
               e_2(t_h) \\
               \vdots \\
               e_{N^{\mathrm{dof}}_{\mathrm{int}}}(t_h)
            \end{bmatrix}_{N^{\mathrm{dof}}_{\mathrm{int}} \times 1},
            \bold{R}_h:=
            \begin{bmatrix}
               r_1(t_h) \\
               r_2(t_h) \\
               \vdots \\
               r_{N^{\mathrm{dof}}_{\mathrm{int}}}(t_h)
            \end{bmatrix}_{N^{\mathrm{dof}}_{\mathrm{int}} \times 1},
        \end{equation*}
        where $r_j(t)=(r(t))_j$.
        Then, following \cite{lui2021spectral,lui2020chebyshev,lui2017legendre}, we deduce that
        \begin{equation*}
            \bA_{Lap}\, \Theta_h= \bold{R}_h, 
        \end{equation*}
        with $|\bold{R}_h|_{\infty} \leq C N^{4} e^{-C N}$ for some positive constant $C$, where $\bA_{Lap}$ is given in \eqref{eq:wave_final}.
        Hence, by \cite[Theorem~3.9]{lui2020chebyshev} and the bounds on $\bA_{Lap}$ from the
        previous subsection, we obtain the desired error estimate.
    \end{proof}

Next, we focus on the convergence of the magnetic field variable $\bB$.

\begin{theorem}\label{Th:B}
    Let $(\bE,\bB)$ be the solution of the Maxwell system
    \eqref{model:original:1}--\eqref{model:original:4} (equivalently, let $\bE$ solve the vector
    wave equation \eqref{hyperbolc_2} and $\bB$ be defined by \eqref{B:decoupled}). Assume that
    both the electric field $\bE$ and the magnetic field $\bB$ are analytic functions in
    $t,x,y$, and $z$. Let $\Gamma_h$ be the error associated with the collocation approximation
    of $\bB$. Then
    \begin{equation}
        \bigl| W^{1/2}\Gamma_h \bigr|_2 \leq c\, N ^{7} e^{-C N}.
    \end{equation}
\end{theorem}

\begin{proof}
        Following the arguments in Theorem~\ref{Th:E} and
        \cite[Theorem~3.9]{lui2020chebyshev}, we derive
        \begin{equation}
            \bD_t \Gamma_h= \widetilde{\bold{R}}_h,
        \end{equation}
        with $|\widetilde{\bold{R}}_h|_{\infty} \leq C N^{6.5} e^{-C N}$,
        where $C$ is a positive constant. Then, using the technique described in
        \cite{lui2021spectral}, we obtain the desired result.
    \end{proof}

The advantage of the proposed method is that we compute the magnetic field as
$\bB_h(t,\bx)=\nabla_h \times \bE_h(t,\bx)$, where $\bE_h(t,\bx)$ is the discrete electric field.
This approach does not enforce the divergence-free constraint of the electric field $\bE_h$
explicitly. However, we have the discrete identity $\nabla_h \cdot (\nabla_h \times \cdot) \approx 0$
and, consequently, the divergence-free condition for the magnetic field is satisfied
asymptotically. Furthermore, we also highlight that the Gauss law for the electric field,
$\|\nabla \cdot \bE_h(t,\bx)- \rho \|_{0,\Omega}$, converges to zero asymptotically.

The next theorem formalizes these observations by showing that the quantities
$\nabla \cdot \bE - \rho$ and $\nabla \cdot \bB$ are not exactly zero but converge to zero with
exponential rate.

\begin{theorem} \label{div_B_E}
    Let $(\bE,\bB)$ be the solution of the Maxwell system
    \eqref{model:original:1}--\eqref{model:original:4}. Assume that both the electric field $\bE$
    and the magnetic field $\bB$ are analytic functions in $t,x,y$, and $z$. Further, assume that
    $\bE_h$ and $\bB_h$ are the solutions of \eqref{hyperbolc_discrete} and
    \eqref{B:decoupled:decoupled}, respectively. Then the following estimates hold:
    \begin{equation*}
        \| \nabla \cdot \bE_h-\rho_h \|_{0,\Omega} \leq c\, N^{6.5} e^{-C N},
    \end{equation*}
    and
    \begin{equation*}
        \| \nabla \cdot \bB_h \|_{0,\Omega} \leq c\, N^9 e^{-C N}.
    \end{equation*}
\end{theorem}

\begin{proof}
    A straightforward calculation yields
    \begin{equation}
    \begin{split}
        \| \nabla \cdot \bE_h-\rho_h \|_{0,\Omega}
        & \leq \| \nabla \cdot \bE_h- \nabla \cdot \bE\|_{0,\Omega}
              + \|\rho-\rho_h \|_{0,\Omega} \\
        & \leq \|\bE_h-\bE \|_{1,\Omega}+\| \rho-\rho_h \|_{0,\Omega} \\
        & \leq c N^2~\|\bE_h-\bE \|_{0,\Omega}+\| \rho-\rho_h \|_{0,\Omega} \\
        & \leq c N^{6.5} e^{-C N}+c e^{-C N} \\
        & \leq c N^{6.5} e^{-C N}.
        \end{split}
    \end{equation}
    Following the same arguments, we rewrite the estimate for the magnetic field as 
    \begin{equation}
    \begin{split}
        \| \nabla \cdot \bB_h \|_{0,\Omega}
        &\leq \| \nabla \cdot \bB_h - \nabla \cdot \bB \|_{0,\Omega} \\
        &\leq C \| \bB_h -  \bB \|_{1,\Omega} \\
        &\leq C N^2~ \| \bB_h -  \bB \|_{0,\Omega} \leq C N^{2} N^{7} e^{-C N}
         = C N^{9}e^{-C N}, 
        \end{split}
    \end{equation}
    where we have used Theorem~\ref{Th:E} and Theorem~\ref{Th:B}, together with the inverse
    estimate provided in \cite[Section~2]{lui2021spectral}.
\end{proof}

In the next section, we describe the tensor-train modification of the full-grid scheme.
\section{Tensor networks}
\label{sec:Tensor-Networks}
In this section, we introduce the \TTf{} format, the specific tensor
network we use in this work, the representation of linear operators in
TT-matrix format, and the cross-interpolation method.
All these methods are fundamental in the tensorization of our spectral
collocation discretization of the Maxwell equations.
For a more comprehensive understanding of notation and concepts, we
refer the reader to the following references:
Refs.~\cite{dolgov2021guaranteed,kolda2009tensor,oseledets2011tensor}, which
provide detailed explanations.

\subsection{Tensor train representation}
\label{subsec:TNs}
\begin{figure}[http] 
  \centering
  \includegraphics[width=0.75\textwidth]{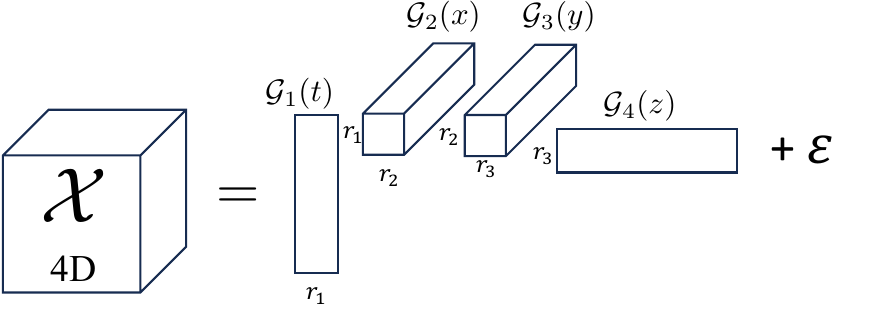}
  \caption{\TTf{} decomposition of a 4D tensor $\ten{X}$, with TT-rank
    $\mathbf{r} = \big[r_1,r_2,r_3\big]$ and approximation error
    $\varepsilon$, in accordance with Eq.~\eqref{eqn:TT_def_element}.}
   \label{fig:TT_4D}
\end{figure}

The TT format, introduced by Oseledets in
2011~\cite{oseledets2011tensor}, represents a sequential chain of
matrix products involving both two-dimensional matrices and
three-dimensional tensors, referred to as TT-cores.
We can visualize this chain as in Figure~\ref{fig:TT_4D}.
Given that tensors in our formulation are at most four-dimensional
(one temporal and three spatial dimensions), we consider the tensor
train format in the context of 4D tensors.
Specifically, the \TTf{} approximation $\ten{X}^{TT}$ of a
four-dimensional tensor $\ten{X}$ is a tensor with elements
\begin{align}
  \ten{X}^{TT}(i_1,i_2,i_3,i_4)
  &=  \sum^{r_1}_{\alpha_{1}=1}\sum^{r_2}_{\alpha_{2}=1}\sum^{r_3}_{\alpha_{3}=1}
  \ten{G}_1(1,i_1,\alpha_1)\ten{G}_2(\alpha_1,i_2,\alpha_2)\ten{G}_3(\alpha_2,i_3,\alpha_3)\ten{G}_4(\alpha_{3},i_4,1)
  + \varepsilon,
  \label{eqn:TT_def_element}
\end{align}
where the error $\varepsilon$ is a tensor with the same dimensions
as $\ten{X}$, and the elements of the array $\mat{r} = [r_1,r_2,r_3]$
are the TT-ranks, which quantify the compression effectiveness of the TT
approximation.
Since each \TTf{} core $\ten{G}_p$ only depends on a single index
of the full tensor $\ten{X}$ (e.g.\ $i_k$), the \TTf{} format effectively
embodies a discrete separation of variables \cite{bachmayr2016tensor}.
In Figure~\ref{fig:TT_4D} we show a four-dimensional array
$\ten{X}(t,x,y,z)$ decomposed in TT-format.

\subsection{Linear operators in TT-matrix format}
\label{SUB:Lin_operinTT}
Suppose that the approximate electric field solution of
\eqref{hyperbolc_compon_1} is represented as a 4D tensor $\ten{V}$.
Then the linear operator $\ten{A}$ acting on that solution is
represented as an 8D tensor.
The transformation $\ten{A}\ten{V}$ is defined as:
\begin{align*}
  \big(\ten{A}\ten{V}\big)(i_1,i_2,i_3,i_4)
  = \sum_{j_1,j_2,j_3,j_4} \ten{A}(i_1,j_1,\dots,i_4,j_4)\,\ten{V}(j_1,\dots,j_4).
\end{align*}
The tensor $\ten{A}$ and the matrix operator $\bA_{Lap}$, defined in
Eq.~\eqref{eq:wave_final}, are related as:
\begin{equation}
  \ten{A}(i_1,j_1,\dots,i_4,j_4)
  = \bA_{Lap}(i_1i_2i_3i_4,j_1j_2j_3j_4).
  \label{eqn:ten_mat_op_relation}
\end{equation}
Thus, we can construct the tensor $\ten{A}$ by suitably reshaping and
permuting the dimensions of the matrix $\bA_{Lap}$.
The linear operator $\ten{A}$ can be further represented in a
variant of \TTf{} format, called \emph{TT-matrix},
cf.~\cite{truong2023tensor}.
The component-wise TT-matrix $\ten{A}^{TT}$ is defined as:
\begin{equation}
  \ten{A}^{TT}(i_1,j_1,\ldots,i_4,j_4)
  =  \sum_{\alpha_{1},\alpha_2,\alpha_{3}}
  \ten{G}_1\big(1,(i_1,j_1),\alpha_1\big)
   \ldots\ten{G}_4\big(\alpha_{3},(i_4,j_4),1\big),
  \label{eqn:TT-matrix-componentwise}
\end{equation}
where $\ten{G}_{k}$ are 4D TT-cores.
\begin{figure}[http]
  \centering
  \includegraphics[width=0.8\textwidth]{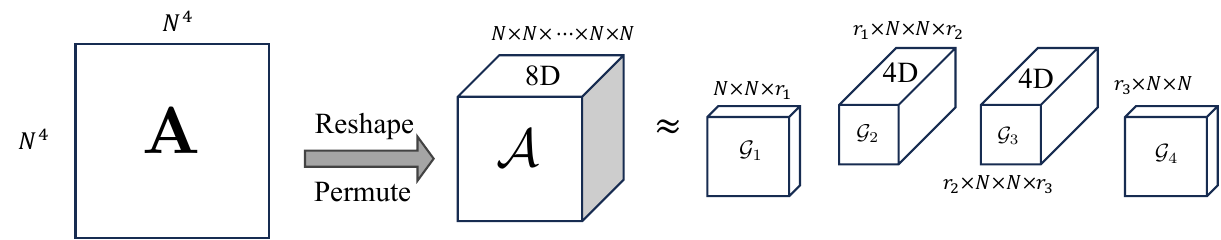}
  \caption{Representation of a linear matrix $\bA_{Lap}$ in the
    TT-matrix format. First, we reshape the operator matrix
    $\bA_{Lap}$ and permute its indices to create the tensor
    $\ten{A}$. Then, we factorize the tensor in the tensor-train
    matrix format according to Eq.~\eqref{eqn:TT-matrix-componentwise}
    to obtain $\ten{A}^{\TTf{}}$.}
  \label{fig:TT-matrix-format}
\end{figure}
Figure~\ref{fig:TT-matrix-format} shows the process of transforming a
matrix operator $\bA_{Lap}$ to its tensor format $\ten{A}$ and,
finally, to its TT-matrix format $\ten{A}^{\TTf{}}$.

We can further simplify the TT-matrix representations of the matrix
$\bA_{Lap}$ if it is a Kronecker product of matrices, i.e.
$\bA_{Lap}=\mat{A}_1\otimes\mat{A}_2\otimes\mat{A}_3\otimes\mat{A}_4$.
Based on the relationship defined in
Eq.~\eqref{eqn:ten_mat_op_relation}, the tensor $\ten{A}$ can be
constructed using tensor products as
$\ten{A} = \mat{A}_1 \tenprod
\mat{A}_2 \tenprod \mat{A}_3 \tenprod \mat{A}_4$.
This implies that the internal ranks of the TT-format of $\ten{A}$
in \eqref{eqn:TT-matrix-componentwise} are all equal to $1$.
In such a case, all summations in
Eq.~\eqref{eqn:TT-matrix-componentwise} reduce to a sequence of single
matrix--matrix multiplications, and the \TTf{} format of $\ten{A}$ becomes
the tensor product of $d$ matrices:
\begin{align}
  \ten{A}^{TT} &= \mat{A}_1 \tenprod \mat{A}_2 \tenprod \dots \tenprod \mat{A}_d.
  \label{eqn:TT_matrices}
\end{align}
This specific structure appears quite often in the matrix
discretization and will be exploited in the tensorization to
construct efficient \TTf{} representations.

\subsection{TT cross interpolation}
\label{sec:CrossInterpolation}
The original \TTf{} algorithm is based on consecutive applications of
singular value decompositions (SVDs) on unfoldings of a tensor
\cite{oseledets2011tensor}.
Although known for its efficiency, the \TTf{} algorithm requires access to the
full tensor, which is impractical and even impossible for extra-large
tensors.
To address this challenge, the cross interpolation algorithm,
TT-cross, has been developed \cite{oseledets2010tt}.
The idea behind TT-cross is essentially to replace the SVD in the TT
algorithm with an \emph{approximate} version of the skeleton/CUR
decomposition \cite{goreinov1997theory, mahoney2009cur}.
CUR decomposition approximates a matrix by selecting a few of its
columns $\mat{C}$, a few of its rows $\mat{R}$, and a matrix $\mat{U}$
that connects them.
Mathematically, CUR decomposition finds an approximation for a matrix
$\mat{A}$ as $\mat{A}\approx \mat{C}\mat{U}\mat{R}$.
The TT-cross algorithm utilizes the maximum volume principle
(\emph{maxvol algorithm})
\cite{goreinov2010find,mikhalev2018rectangular} to determine
$\mat{U}$.
The maxvol algorithm chooses a few columns $\mat{C}$ and rows
$\mat{R}$ of $\mat{A}$ such that the intersection matrix
$\mat{U}^{-1}$ has maximum volume \cite{savostyanov2011fast}.

TT-cross interpolation and its variants can be seen as heuristic
generalizations of CUR to tensors
\cite{oseledets2008tucker,sozykin2022ttopt}.
TT-cross utilizes the maximum volume algorithm iteratively, often
beginning with a few randomly chosen fibers, to select an optimal set
of specific tensor fibers that capture the essential information of the
tensor \cite{savostyanov2014quasioptimality}.
These fibers are used to construct a lower-rank \TTf{}
representation. The naive generalization of CUR is known to be
expensive, which has led to the development of various heuristic
optimization techniques such as TT-ALS \cite{holtz2012alternating},
DMRG \cite{oseledets2012solution,savostyanov2011fast}, and AMEN
\cite{dolgov2014alternating}.
To solve the Maxwell equation in space--time format, we use TT-cross to build the \TTf{} format
directly from the boundary conditions, initial
conditions, and loading functions.

\subsection{Tensorization of the electric field}
\label{tensor:electric}
From Section~\ref{sec:wave_equation}, we can immediately identify the
core tensors for the TT-modification of the space--time discrete system.
The space--time discretization produces a linear system for all
interior nodes as specified in Eq.~\eqref{eq:wave_final}.
Here, we simplify the notation and refer to this equation in the compact form
$\bA_{Lap} \boldsymbol{V}=\boldsymbol{F} -\boldsymbol{F}^{\text{BD}}$, where
$\boldsymbol{V}=\hat{\bv}_h$ collects the interior unknowns corresponding to the electric field component
under consideration.
Tensorization is the process of building the \TTf{} format of all components
of this linear system,
\begin{equation}
  \bA_{Lap} \boldsymbol{V}
  = \boldsymbol{F} -\boldsymbol{F}^{\text{BD}}
  \;\;\longrightarrow\;\;
  \ten{A}^{TT}\ten{V}^{TT}
  = \ten{F}^{TT} - \ten{F}^{\text{BD},TT},
\end{equation}
where $\ten{A}^{TT} = \ten{A}_t^{TT} + \ten{A}_D^{TT}$.
In matrix form, the operator $\bA_{Lap}$ acts on the vectorized solution
$\boldsymbol{V}$.
In the full tensor format, the solution is kept in its original tensor
form $\ten{V}$, which is a 4D tensor.
Consequently, the operators $\ten{A}_t$ and $\ten{A}_D$ are 8D tensors.
Lastly, $\ten{F}$ and $\ten{F}^{\text{BD}}$ are both 4D tensors.
Given that these operators in matrix form have Kronecker
product structure, their \TTf{} format can be constructed by using
component matrices as \TTf{} cores.
To construct the tensor format of the operators acting on the interior
nodes, we define the index sets
\begin{equation*}
  \begin{split}
    & \ten{I}_t = 2:N \quad \text{(index set for the time variable)},\\
    & \ten{I}_s = 2:(N-1) \quad \text{(index set for each spatial variable)}.
  \end{split}
\end{equation*}

\PGRAPH{$\bullet$ TT-matrix time operator} $\ten{A}^{TT}_t$: The temporal operator
in TT-matrix format acting only on the interior nodes is constructed
as:
\begin{equation}
  \ten{A}^{TT}_t
  = \bI_{N-1} \tenprod \bI_{N-1} \tenprod \bI_{N-1} \tenprod \bS_t(\ten{I}_t,\ten{I}_t),
  \label{eqn:At_TT_op}
\end{equation}
where $\bI_{N-1}$ is the identity matrix of size
$(N-1)\times(N-1)$.
  
\PGRAPH{$\bullet$ TT-matrix diffusion operator} $\ten{A}^{TT}_D$: The Laplace
operator in TT-matrix format is constructed as:
\begin{equation}
  \begin{split}
      \ten{A}^{TT}_D
      &= \bS_{xx}(\ten{I}_s,\ten{I}_s) \tenprod \bI_{N-1}\tenprod \bI_{N-1} \tenprod \bI_{N}
      + \bI_{N-1} \tenprod \bS_{yy}(\ten{I}_s,\ten{I}_s)  \tenprod\bI_{N-1} \tenprod \bI_{N} \\
      &\quad + \bI_{N-1} \tenprod  \bI_{N-1} \tenprod\bS_{zz}(\ten{I}_s,\ten{I}_s) \tenprod \bI_N .
  \end{split}
\end{equation}

\PGRAPH{$\bullet$ TT loading tensor} $\ten{F}^{TT}$: The TT loading tensor
$\ten{F}^{TT}$ is constructed by applying TT-cross interpolation
to the function $f(t,x,y,z)$ on the grid of interior
nodes. Such techniques are investigated in
\cite{adak2024tensorA,adak2024tensor}. We highlight that if the forcing
term satisfies a \emph{separation-of-variables} condition, cross interpolation
may not be necessary to construct the loading tensor.

\PGRAPH{$\bullet$ TT boundary tensor}  $\ten{F}^{\text{BD},TT}$: Unlike full-grid computation,
the implementation of boundary conditions for a tensor-train grid is not straightforward.
In this direction, we construct an operator  $\ten{A}^{\text{map}}$ which imposes the given boundary
and initial data at the correct positions. In particular, $\ten{A}^{\text{map}}$ maps all nodes
to the subset of unknown nodes.
The size of the operator $\ten{A}^{\text{map}}$ is
$(N-2)\times (N) \times
(N-1)\times (N+1) \times (N-1)\times (N+1) \times(N)\times
(N+1)$ according to the . Here, odd indices signify only the number of unknown nodes,
while even indices signify the total number of nodes. 
The boundary tensor $\ten{F}^{\text{BD}}$ is then computed by applying
$\ten{A}^{\text{map}}$ to a tensor $\ten{G}^{\text{BD}}$ containing only the
information from boundary and initial conditions.
The details about constructing $\ten{A}^{\text{map},TT}$ and
$\ten{G}^{\text{BD},TT}$ are included in \ref{APP:A_map_tt}. Here, we focused on the tensorization of $E_x$ variable corresponding to \eqref{hyperbolc_compon_1}. The tensorization for other components such as $E_y$, and $E_z$ will be followed same technique as mentioned above.

At this point, we have completed the construction of the TT-format of the
linear system $\ten{A}^{TT}\ten{V}^{TT} = \ten{F}^{TT} -
\ten{F}^{\text{BD},TT}$.
To solve this \TTf{} linear system by optimization techniques, we use
the routines \texttt{amen\_cross}, \texttt{amen\_solve} and \texttt{amen\_mm}
from the MATLAB TT-Toolbox \cite{Oseledets:2023}.

\subsection{Tensorization of the magnetic field}
\label{tensor:magnetic:field}
In this section, we primarily focus on the TT format of the magnetic-field component
$B_{xh}$ in \eqref{B:decoupled:decoupled}; the computation of the other components
$B_{yh}$ and $B_{zh}$ is analogous. The TT-modification of the interior system
\eqref{eq:magnetic_interior} can be written as
\begin{equation}
    \ten{A}^{TT}_t  \ten{B}_x^{TT}
    = \mathcal{F}_{\text{curl}}^{TT} - \mathcal{F}_{\text{curl}}^{\text{BD},TT},
\end{equation}
where $\ten{A}^{TT}_t$ is defined in \eqref{eqn:At_TT_op}, and $\ten{B}_x^{TT}$ is the 4D tensor
representing the magnetic-field component $B_{xh}$. Further, $\mathcal{F}_{\text{curl}}^{TT}$ is a
4D tensor computed as
\begin{equation}
\begin{split}
  \mathcal{F}_{\text{curl}}^{TT}
  &=  \Big ( \bI_{x}(\ten{I}_s,\ten{I}_s) \tenprod \bS_{y}(\ten{I}_s,\ten{I}_s) \tenprod \bI_{z}(\ten{I}_s,\ten{I}_s) \tenprod \bI_{t}(\ten{I}_t,\ten{I}_t) \Big ) \boldsymbol{\mathcal{I}}_{V_z \to S_x} \boldsymbol{E}_z \\
  & \quad - \Big ( \bI_{x}(\ten{I}_s,\ten{I}_s) \tenprod  \bI_{y}(\ten{I}_s,\ten{I}_s)  \tenprod \bS_{z}(\ten{I}_s,\ten{I}_s) \tenprod \bI_{t}(\ten{I}_t,\ten{I}_t) \Big ) \boldsymbol{\mathcal{I}}_{V_y \to S_x}   \boldsymbol{E}_y,
  \end{split}
\end{equation}
where $\boldsymbol{\mathcal{I}}_{V_z \to S_x}$ is the matrix corresponding to the interpolation operator
that maps the vector of collocation values $\boldsymbol{E}_z$ to the grid associated with $V_k^x$
in \eqref{grid:B}. Similarly, $\boldsymbol{\mathcal{I}}_{V_y \to S_x}$ maps the vector $\boldsymbol{E}_y$
to the grid corresponding to $V_k^x$. We highlight that this interpolation operator is needed since
we have computed the electric and magnetic field components on different staggered grids. In this
post-processing technique, one must also treat boundary nodes carefully; this can be done following
the techniques in \cite{adak2024tensorA}.

Further, we need to tensorize the boundary term $\mathcal{F}_{\text{curl}}^{\text{BD},TT}$, which is given by
\begin{equation}
    \mathcal{F}_{\text{curl}}^{\text{BD},TT}
    = \bI_{x}(\ten{I}_s,:) \tenprod \bI_{y}(\ten{I}_s,:) \tenprod \bI_{z}(\ten{I}_s,:) \tenprod \bS_{t}(\ten{I}_t,:) \,\mathcal{G}_{B_x}^{\text{BD},TT}, 
\end{equation}
where $\mathcal{G}_{B_x}^{\text{BD},TT}$ is the $(N+1) \times N \times N \times (N+1)$ boundary tensor
that consists only of boundary and initial values and zeros at the interior grid points. For a more
detailed discussion, we refer to \cite{adak2024tensor,oseledets2011tensor} and
\ref{APP:A_map_tt}.
\section{Numerical Experiments}
\label{sec:numerical_experiments}

In this section, we present numerical experiments that assess the performance of our tensor-network–based space–time spectral collocation method for Maxwell's equations in Tensor-Train (TT) format. The first two tests use manufactured solutions of rank~1, while the third test considers a rank~3 solution. In all experiments, the tolerance used in the TT truncation is chosen to be smaller than the approximation error introduced when constructing the TT representation. Throughout, we compare the TT-based solver with a full-grid solver that operates directly on dense tensors. All simulations are performed in \textsc{Matlab} on a Mac system with an M2 processor.

\subsection{Experiment~1} \label{Ex:1}

In the first test, we apply the TT-solver to \eqref{model:original:1}–\eqref{model:original:4} using the manufactured solution
\begin{equation}
\label{Ex1}
    \bE(\bx,t)=\begin{bmatrix}
               0 \\
               \sin(2\pi x)\,\sin(2\pi y)\,\sin(2\pi t) \\
               \sin(\pi x)\,\sin(\pi y)\,\cos(\pi t)
            \end{bmatrix},
    \qquad
    \bB(\bx,t)=\begin{bmatrix}
               -\sin(\pi x)\,\cos(\pi y)\,\sin(\pi t) \\
                \cos(\pi x)\,\sin(\pi y)\,\sin(\pi t) \\
                \cos(2\pi x)\,\sin(2\pi y)\,\cos(2\pi t)
            \end{bmatrix}.
\end{equation}
The computational domain is $\Omega_T=[0,1]\times[0,1]^3$, with $\epsilon_0 = 1$ and $c = 1$. The exact solution is independent of the $z$-variable, and boundary and initial data are chosen accordingly.

\begin{figure}[htbp]
  \begin{center}
    \begin{tabular}{cc}
      \includegraphics[width=0.4\linewidth]{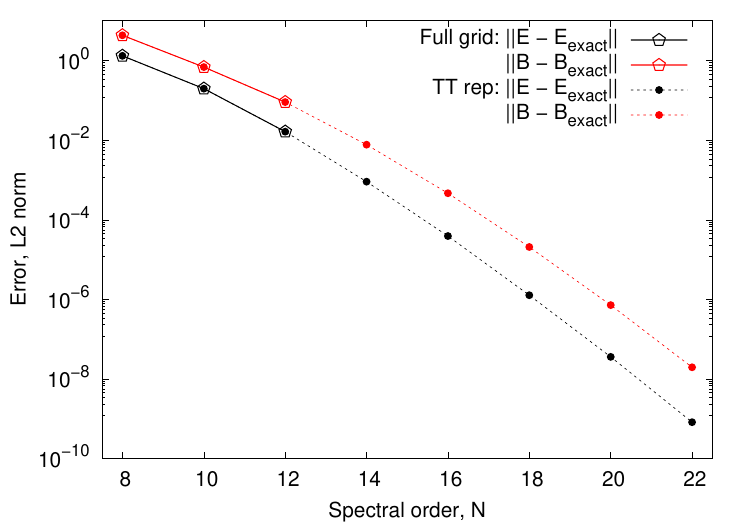}
      \includegraphics[width=0.4\linewidth]{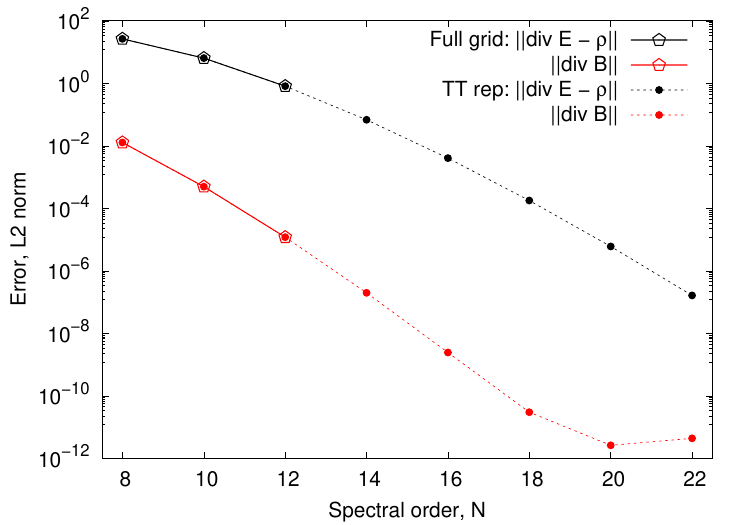}
      \\
      \includegraphics[width=0.4\linewidth]{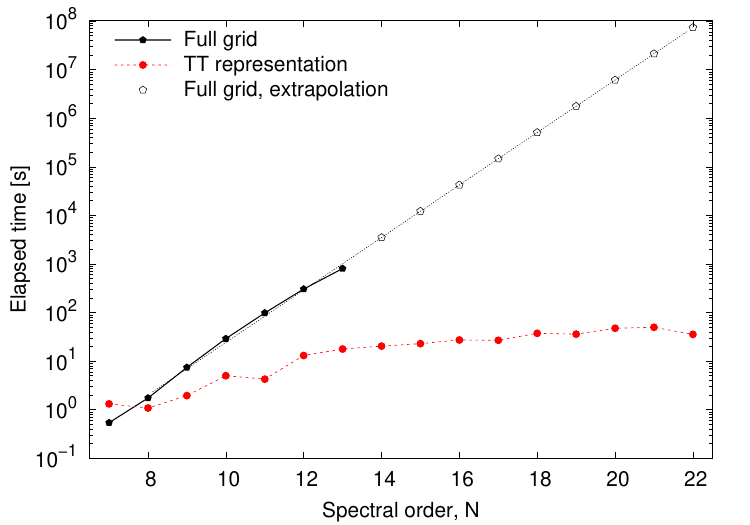}
    \end{tabular}
    \caption{Test Case~1, upper panels: $L_2$-errors in the computation of the electric and magnetic fields (upper left), and divergence constraints: $\|\div \bE - \rho\|$, $\| \div \bB \|$ (upper right).
    Bottom panel: elapsed time for the full grid calculation vs TT representation. 
    }
    \label{Ex1:con:B_E}
  \end{center}
\end{figure}

The TT-solver matches the accuracy of the full-grid solver, reaching errors on the order of
$10^{-10}$ with 22 grid points per dimension when using a TT tolerance of $10^{-11}$. The
convergence behaviour is illustrated in Figure~\ref{Ex1:con:B_E}: the upper-left panel shows the
$L_2$-errors in the electric and magnetic fields as a function of the number of grid points,
while the upper-right panel shows the corresponding divergence errors
$\|\div \bE - \rho\|$ and $\|\div \bB\|$. The divergence of the magnetic field converges
exponentially to zero, consistent with the fact that $\bB$ is computed in the curl space of
$\bE$ and $\div(\curl \cdot)=0$, and Gauss’s law for $\bE$ is also satisfied with exponential
accuracy. These observations confirm the theoretical results of Theorem~\ref{div_B_E}.

Finally, we highlight the difference in computational cost. The bottom panel of
Figure~\ref{Ex1:con:B_E} compares the elapsed time of the TT-based and full-grid solvers. While
both solvers initially require comparable time, the full-grid solver becomes significantly more
expensive as the grid is refined; for example, at 22 points per dimension the TT-solver is
approximately $10^7$ times faster.

\subsection{Experiment~2}

In the second experiment, we again study \eqref{model:original:1}–\eqref{model:original:4}, now with a three–dimensional manufactured solution:
\begin{equation}
    \bE(\bx,t)=\begin{bmatrix}
               0 \\
               \sin(2\pi x)\,\sin(2\pi y)\,\sin(2\pi z)\,\sin(2\pi t) \\
               \sin(\pi x)\,\sin(\pi y)\,\sin(\pi z)\,\cos(\pi t)
            \end{bmatrix},
\end{equation}
\begin{equation}
    \bB(\bx,t)=\begin{bmatrix}
               -\sin(\pi x)\,\cos(\pi y)\,\sin(\pi z)\,\sin(\pi t)
               -\sin(2\pi x)\,\sin(2\pi y)\,\cos(2\pi z)\,\cos(2\pi t) \\
               \cos(\pi x)\,\sin(\pi y)\,\sin(\pi z)\,\sin(\pi t) \\
               \cos(2\pi x)\,\sin(2\pi y)\,\sin(2\pi z)\,\cos(2\pi t)
            \end{bmatrix}.
\end{equation}
The domain is $\Omega_T=[0,1]\times[-1,1]^3$ with $\epsilon_0=1$ and $c=1$. Initial and boundary data are derived from the exact solution.

\begin{figure}[htbp]
  \begin{center}
    \begin{tabular}{cc}
      \includegraphics[width=0.4\linewidth]{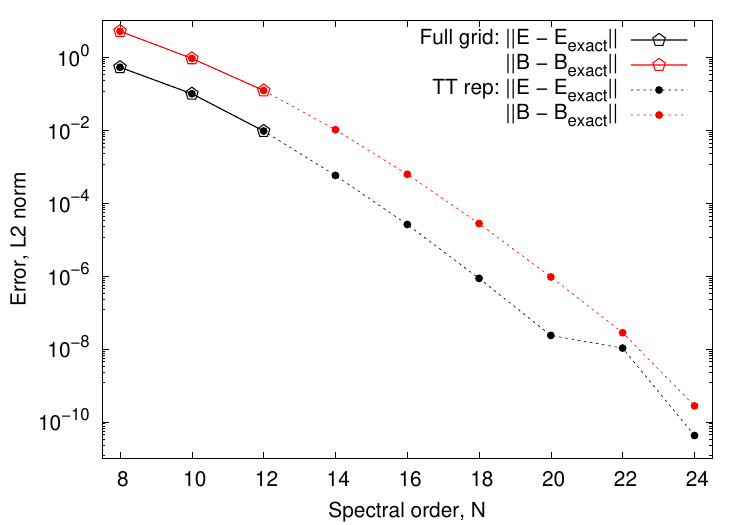}
      \includegraphics[width=0.4\linewidth]{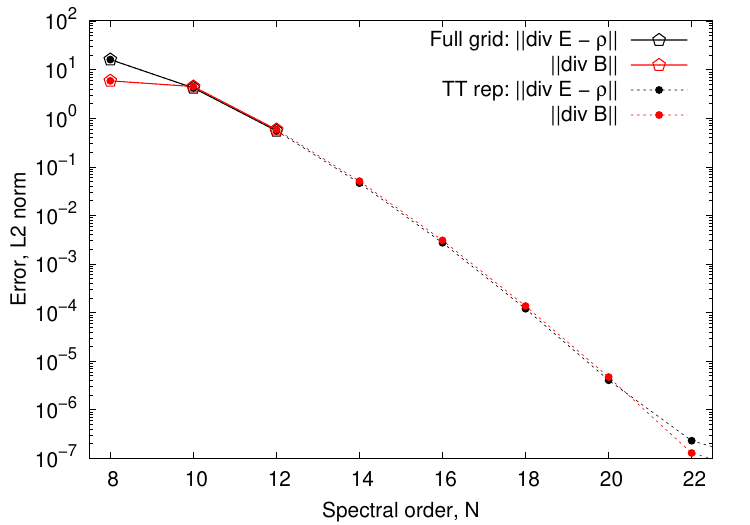}
      \\
      \includegraphics[width=0.4\linewidth]{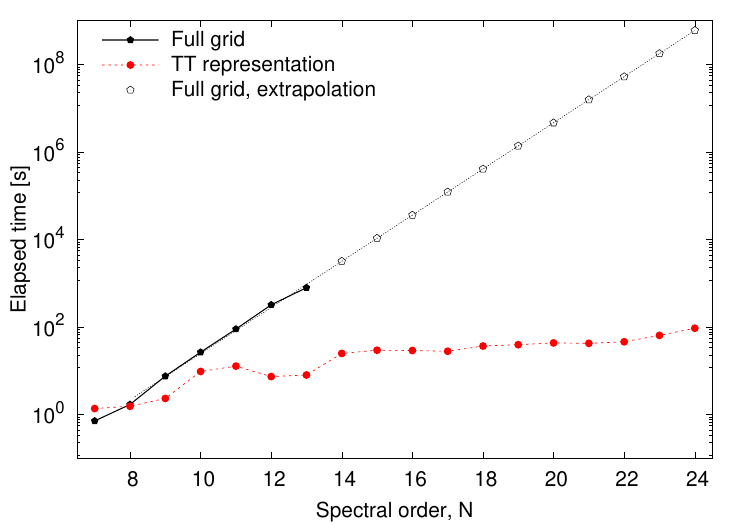}
    \end{tabular}
    \caption{Test Case~2, upper panels: $L_2$-errors of electric and magnetic fields (upper left),
    and divergence constraints: $\|\div \bE - \rho \|$, $\| \div \bB \| $ (upper right).
    Bottom panel: elapsed time for full grid calculation vs TT representation.}
    \label{Ex2:con:B_E}
  \end{center}
\end{figure}

The TT-solver again reproduces the accuracy of the full-grid solver on coarse and moderate grids
(see Figure~\ref{Ex2:con:B_E}). Beyond 12 points per dimension, the full-grid method exceeds
available memory, whereas the TT-solver remains efficient and stable up to at least 22 points
per dimension. The divergence of $\bB$ converges exponentially to zero, and Gauss’s law for $\bE$
is satisfied with the same rate, again matching the theoretical predictions of
Theorem~\ref{div_B_E}.

Timing results in Figure~\ref{Ex2:con:B_E} show the full-grid solver scaling steeply with grid size,
while the TT-based solver exhibits a much milder growth in cost. Extrapolation of full-grid timings
indicates that for a 24-point grid, the TT-solver is roughly $10^8$ times faster.

\subsection{Experiment~3}

The third test considers a higher-rank solution:
\begin{equation}
    \bE(\bx,t)=\begin{bmatrix}
               0 \\
               0 \\
               \displaystyle\sum_{k=1}^{3}\sin(k\pi x)\,\sin(k\pi y)\,\cos(k\pi t)
            \end{bmatrix},
    \qquad
    \bB(\bx,t)=\begin{bmatrix}
               -\displaystyle\sum_{k=1}^{3}\sin(k\pi x)\,\cos(k\pi y)\,\sin(k\pi t) \\
                \displaystyle\sum_{k=1}^{3}\cos(k\pi x)\,\sin(k\pi y)\,\sin(k\pi t) \\
               0
            \end{bmatrix}.
\end{equation}
Boundary and initial conditions are taken from the exact expressions. This experiment evaluates performance when the true solution has higher tensor rank.

\begin{figure}[htbp]
  \begin{center}
    \begin{tabular}{cc}
      \includegraphics[width=0.4\linewidth]{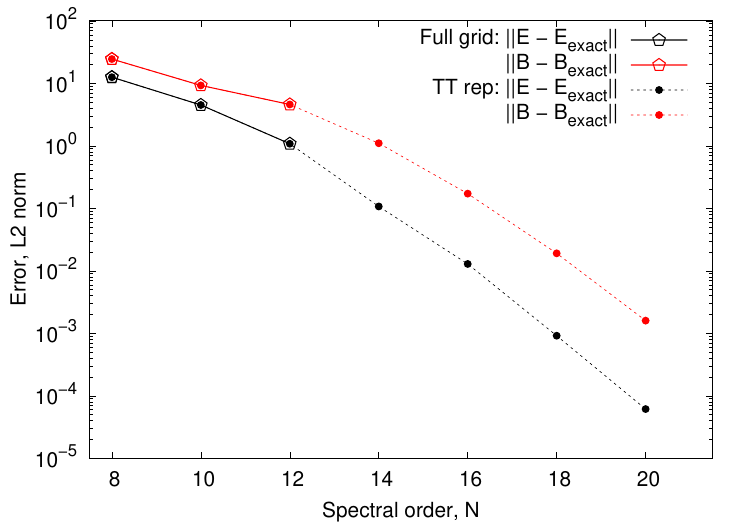}
      \includegraphics[width=0.4\linewidth]{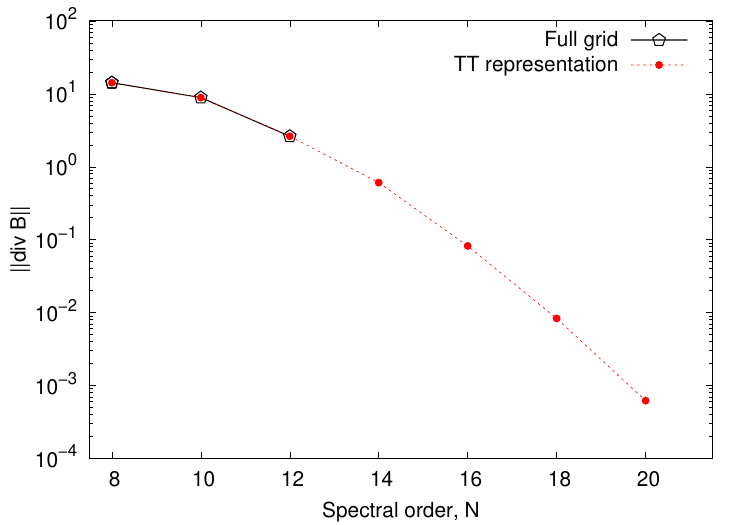}
      \\
      \includegraphics[width=0.4\linewidth]{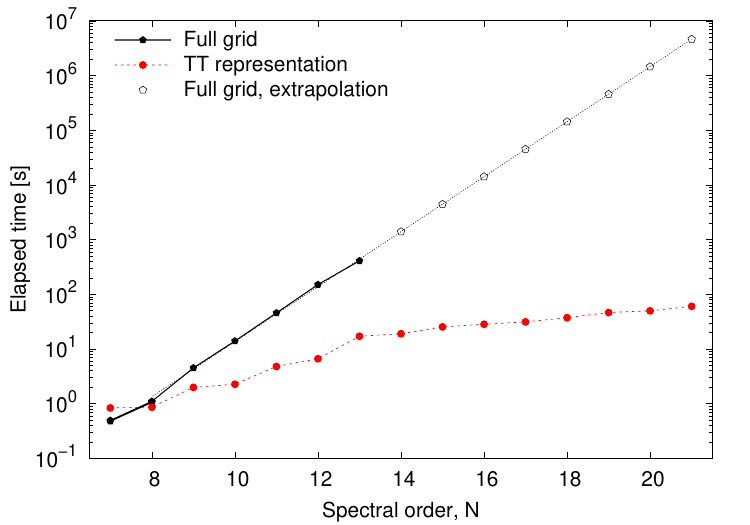}
    \end{tabular}
    \caption{Test Case~3, upper panels: $L_2$-error in the computation of the electric and
    magnetic fields (upper left), and divergence constraint $\div \bB = 0$ (upper right).
    Bottom panel: elapsed time for full-grid calculation vs.\ TT representation.}
    \label{Ex3:con:B_E}
  \end{center}
\end{figure}

The TT- and full-grid solvers deliver comparable accuracy across all tested grid sizes
(see Figure~\ref{Ex3:con:B_E}). As in the previous experiments, both methods require similar
computational effort on coarse grids, but the full-grid cost grows rapidly with refinement.
At 22 grid points per dimension, the TT-solver is about $10^7$ times faster than the full-grid
method. These results confirm that the TT approach remains effective even for solutions with
higher intrinsic rank, and that it is well suited for large-scale Maxwell simulations involving
many variables.

\section{Conclusions}
\label{sec:conclusions}
In this work, we have developed and analyzed a space-time spectral collocation tensor network method for the efficient numerical solution of Maxwell's equations in three spatial dimensions. Our approach combines the exponential convergence properties of Chebyshev spectral methods with the dimensional scalability of tensor-train decomposition, addressing the curse of dimensionality that has long hindered high-resolution electromagnetic simulations.

The key contributions of this work are threefold. First, we have formulated a space-time spectral collocation scheme that discretizes both spatial and temporal dimensions using Chebyshev polynomials, yielding a global system that captures the entire time evolution simultaneously. By reformulating Maxwell's equations as a second-order wave equation for the electric field under constant-coefficient assumptions, we decouple the system into independent componentwise problems with natural Kronecker product structure. Second, we have introduced a staggered spectral discretization where electric and magnetic field components live in different polynomial spaces with carefully chosen degrees of freedom in each coordinate direction. This spectral staggering preserves the divergence-free property of the magnetic field to spectral accuracy without explicit constraint enforcement, analogous to classical FDTD staggering but within a high-order spectral framework. Third, we have developed a complete tensorization strategy that exploits the Kronecker structure of the discrete operators to construct efficient tensor-train representations using TT-cross interpolation and alternating minimal energy solvers.

Our theoretical analysis establishes rigorous foundations for the method. We have derived condition number estimates showing that the discrete wave operator scales as $\mathcal{O}(N^4)$ while the temporal operator for magnetic field reconstruction scales as $\mathcal{O}(N^2)$, where $N$ is the polynomial degree. More importantly, we have proven exponential convergence of the discrete electric and magnetic fields in the $L^2$ norm, with convergence rates of the form $\mathcal{O}(N^{p} e^{-CN})$ for appropriate powers $p$. Furthermore, we have established that both Gauss's law for the electric field and the divergence-free constraint for the magnetic field are satisfied with exponential accuracy, confirming that the discrete divergence operator applied to the discrete curl vanishes asymptotically.

Numerical experiments validate these theoretical results and demonstrate the practical efficiency of the method. Across three test cases with manufactured solutions of varying tensor rank, we consistently observe exponential convergence matching the full-grid spectral method while achieving dramatic computational savings. For problems with $22$ collocation points per dimension (corresponding to approximately $10^{8}$ degrees of freedom in the space-time grid), the tensor-train solver is roughly $10^7$ times faster than the full-grid approach. Critically, the TT method extends beyond the memory limitations of full-grid computation: while full-grid methods become infeasible beyond $12$--$14$ points per dimension, the TT solver successfully handles $22$--$24$ points, accessing a regime that would be completely intractable with conventional approaches.

Our work demonstrates that space-time spectral collocation methods, when combined with tensor-train compression, provide a powerful tool for solving Maxwell's equations with unprecedented efficiency and accuracy. By achieving exponential convergence while maintaining linear complexity in storage and computation, the method opens new possibilities for high-fidelity electromagnetic simulations that were previously computationally prohibitive. The rigorous mathematical framework, validated through numerical experiments, establishes a solid foundation for further developments in tensor network methods for computational electromagnetics.

\appendix
\section{Construction of \texorpdfstring{$\ten{A}^{map,TT} \text{ and } 
\label{APP:A_map_tt}
\ten{G}^{TT}$}{} \cite{adak2024tensor}} Here we provide details about the construction of
$\ten{A}_{E_x}^{map,BD}$, $\ten{A}_{E_x}^{map,IC}$, $\ten{A}_{E_{x,t}}^{map,BD}$, $\ten{A}_{E_{x,t}}^{map,IC}$ and $\ten{G}^{BD,TT}_{E_x}$, $\ten{G}^{IC,TT}_{E_x}$, $\ten{G}^{BD,TT}_{E_{x,t}}$, $\ten{G}^{IC,TT}_{E_{x,t}}$, where $BD$, and $IC$ stand for boundary condition and initial condition in space and time variables respectively. The sub-index $E_x$, and $E_{x,t}$ signify implementation boundary condition or initial guess  for $E_x$, and $\frac{\partial E_x}{ \partial t}$ respectively.
%
We define
\begin{equation}
\label{funda:Amap}
\begin{split}
  \ten{A}^{map}_{1} & =   \bI^{N+1}(\ten{I}_s,\ten{I}_s )\tenprod \bI^{N+1}(\ten{I}_s,\ten{I}_s)\tenprod \bI^{N+1}(\ten{I}_s,\ten{I}_s) \tenprod \bold{S}_t(\ten{I}_t,\ten{I}_t) ,\\
  \ten{A}^{map}_{2} & = \bI^{N+1}(\ten{I}_s,: )\tenprod \bI^{N+1}(\ten{I}_s,:)\tenprod \bI^{N+1}(\ten{I}_s,:) \tenprod \bold{S}_t(\ten{I}_t,:), \\
  \ten{A}^{map}_{3} & = \bI^{N+1}(\ten{I}_s,: )\tenprod \bI^{N+1}(\ten{I}_s,:)\tenprod \bI^{N+1}(\ten{I}_s,:) \tenprod \bI^{N+1}(\ten{I}_t,:),\\
  \ten{A}^{map}_{4} & = \bold{S}_{xx}(\ten{I}_s,: )\tenprod \bI^{N+1}(\ten{I}_s,:)\tenprod \bI^{N+1}(\ten{I}_s,:) \tenprod \bI^{N+1}(\ten{I}_t,:) \\
  & \quad + \bI^{N+1}(\ten{I}_s,: )\tenprod \bold{S}_{yy}(\ten{I}_s,:)\tenprod \bI^{N+1}(\ten{I}_s,:) \tenprod \bI^{N+1}(\ten{I}_t,:) \\
  & \quad + \bI^{N+1}(\ten{I}_s,: )\tenprod \bI^{N+1}(\ten{I}_s,:)\tenprod \bold{S}_{zz}(\ten{I}_s,:) \tenprod \bI^{N+1}(\ten{I}_t,:).
\end{split}
\end{equation}
By using \eqref{funda:Amap}, we construct the final mapping as follows
\begin{equation}
\begin{split}
    \ten{A}_{E_x}^{map,BD} & = \ten{A}^{map}_{4}, \qquad \qquad 
    \ten{A}_{E_x}^{map,IC}  = \ten{A}^{map}_{1} \ten{A}^{map}_{2}, \\
     \ten{A}_{E_{x,t}}^{map,BD} & = \ten{A}^{map}_{1}  \ten{A}^{map}_{3},
     \quad \; \; \ten{A}_{E_{x,t}}^{map,IC} = \ten{A}^{map}_{2}.
\end{split}
\end{equation}

Next, we show how to construct the $\ten{G}^{BD}$ tensor.
The tensor $\ten{G}^{BD}$ is a $(N)\times (N+1) \times (N+1) \times
(N+1)$, in which only the BD/IC nodes are computed. Other nodes
are zeros. The TT tensor $\ten{G}^{BD}$ is constructed using the cross
interpolation.
Then, the boundary tensor $\ten{F}^{BD,TT}$ is computed as:
\begin{equation}
\label{boundary tesor:E}
  \ten{F}^{BD,TT}= \ten{A}_{E_x}^{map,BD} \ten{G}^{BD,TT}_{E_x}+ \ten{A}_{E_x}^{map,IC} \ten{G}^{IC,TT}_{E_x} + \ten{A}_{E_{x,t}}^{map,BD} \ten{G}^{BD,TT}_{E_{x,t}} +\ten{A}_{E_x}^{map,IC} \ten{G}^{IC,TT}_{E_x}.
\end{equation}

\section*{Data availability statement}

The data that support the findings of this research are available from the corresponding author upon reasonable request.

\section*{Declaration of competing interest}

The authors declare that they have no known competing financial interests or personal relationships that could have appeared to influence the work reported in this paper.

\section*{CRediT authorship contribution statement}
D. Adak, D. P. Truong, R. Chinomona, O. Korobkin, Kim O. Rasmussen, B. S. Alexandrov : Conceptualization,
Methodology, Writing-original draft, Review $\&$ Editing.


\section*{Acknowledgments}
The authors gratefully acknowledge the support of the Laboratory
Directed Research and Development (LDRD) program of Los Alamos
National Laboratory under project number 20240705ER.
Los Alamos National Laboratory is operated by Triad National Security,
LLC, for the National Nuclear Security Administration of
U.S. Department of Energy (Contract No.\ 89233218CNA000001).


\bibliographystyle{plain}
\bibliography{ttrain}


\clearpage

\end{document}